\documentclass[a4paper,12pt]{article}
\usepackage{amssymb}
\usepackage{verbatim}
\usepackage{amsmath}
\usepackage{amscd}
\usepackage{array}
\usepackage[arrow,matrix]{xy}
\usepackage[bookmarksnumbered,colorlinks,linkcolor=black,citecolor=black,urlcolor=black]{hyperref}

\newtheorem{theo}{Theorem}[section]
\newtheorem{prop}[theo]{Proposition}
\newtheorem{lem}[theo]{Lemma}
\newtheorem{cor}[theo]{Corollary}
\newtheorem{defi}[theo]{Definition}
\newtheorem{rem}[theo]{Remark}
\newtheorem{ex}[theo]{Example}
\newtheorem{con}[theo]{Convention}

\def \deg {{\rm{deg}}}
\def \Br {{\rm{Br}}}

\def \si {{\sigma}}
\def \Ga {{\Gamma}}
\def \R {{\mathbb{R}}}
\def \Pic {{\rm {Pic}}}

\def \Gal {{\rm{Gal}}}
\def \Ker {{\rm{Ker}}}
\def \Im {{\rm {Im}}}

\def \Spec {{\rm{Spec}}}
\def \dim {{\rm{dim}}}
\def \Hom {{\rm {Hom}}}
\def \End {{\rm {End}}}
\def \Pic {{\rm {Pic}}}
\def \GL {{\rm {GL}}}

\def \SO {{\rm {SO}}}

\def \Sp {{\rm {Sp}}}

\def \tr {{\rm {tr}}}
\def \Bi {{\rm {Bi}}}

\def\s{{\rm s}}
\def\ov{\overline}

\def \Z {{\mathbb Z}}
\def \Q {{\mathbb Q}}
\def \F {{\mathbb F}}

\def \Id {{\rm Id}}

\def \Mat {{\rm{Mat}}}

\def \rk {{\rm{rk}}}

\def \Ad {{\rm Ad}}
\def \ad {{\rm ad}}
\def \gl {{\mathfrak{gl}}}
\def \J {{\mathfrak J}}

\def\G{{\mathbb G}}
\def\D{{\cal D}}

\def\C{{\mathbb C}}

\def\lra{\longrightarrow}

\def\id{{\rm id}}

\def\H{{\rm H}}
\def\ad{{\rm ad}}

\def\NS{{\rm NS\,}}

\def\O{{\cal O}}

\def\V{{\cal V}}
\def\W{{\cal W}}

\def\si{\sigma}

\def\Ga{\Gamma}

\def\et{{\rm{\acute et}}}

\def\g{{\mathfrak g}}

\newcommand{\bthe}{\begin{theo}}
\newcommand{\ble}{\begin{lem}}
\newcommand{\bpr}{\begin{prop}}
\newcommand{\bco}{\begin{cor}}
\newcommand{\bde}{\begin{defi}}
\newcommand{\ethe}{\end{theo}}
\newcommand{\ele}{\end{lem}}
\newcommand{\epr}{\end{prop}}
\newcommand{\eco}{\end{cor}}
\newcommand{\ede}{\end{defi}}
\newcommand{\brem}{\begin{rem}}
\newcommand{\erem}{\end{rem}}
\newcommand{\bex}{\begin{ex}}
\newcommand{\eex}{\end{ex}}
\newcommand{\bcon}{\begin{con}}
\newcommand{\econ}{\end{con}}

\usepackage[textheight=220mm,textwidth=150mm]{geometry}

\title{Invariant Brauer group of an abelian variety}
\author{M. Orr, A.N. Skorobogatov, D. Valloni and Yu.G. Zarhin}
\date{\today}
\begin{document}
\maketitle

\begin{abstract}
We study a new object that can be attached to an abelian variety or 
a complex torus: the invariant Brauer group, as
recently defined by Yang Cao. Over the field of complex numbers
this is an elementary abelian 2-group with an explicit upper bound on the rank.
We exhibit many cases in which the invariant Brauer group is zero, and construct complex
abelian varieties in every dimension starting with 2, both simple and 
non-simple, with invariant Brauer group of order~2. 
We also address the situation in finite characteristic and over non-closed fields.
\end{abstract}

\section{Introduction}

The Brauer--Grothendieck group $\Br(X)=\H^2_\et(X,\G_m)$ is an important cohomological invariant of 
an algebraic variety $X$ over a field $k$. It has applications in algebraic geometry (the rationality problem 
over algebraically closed fields) and in diophantine geometry (the Brauer--Manin obstruction to the
local-to-global principles over global fields), see \cite{CTS}. If $k$ is not separably closed, a more accessible part
of $\Br(X)$ is the {\em algebraic Brauer group} $\Br_1(X)$ defined as the kernel of the natural map
$\Br(X)\to\Br(X_\s)$, where $X_\s=X\times_k k_\s$ and $k_\s$ is a separable closure of $k$.

If $X$ is acted on by an algebraic $k$-group, it is natural to look for an `equivariant' version of
the Brauer--Grothendieck group that takes into account the symmetries of $X$. 
Such a version was suggested recently by Y.~Cao in the case of an action
$m\colon G\times X\to X$ of a {\em connected} algebraic $k$-group $G$, see \cite[D\'ef. 3.1]{C18}, 
\cite[D\'ef. 1.1 (2)]{YC}. He defined the {\em invariant Brauer group} $\Br_G(X)$ 
as the subgroup of $\Br(X)$ consisting of the elements $x\in\Br(X)$
such that $m^*(x)\in \Br(G\times X)$ belongs to the subgroup $\pi_1^*\Br(G)+\pi_2^*\Br(X)$,
where $\pi_1\colon G\times X\to G$ and $\pi_2\colon G\times X\to X$ are the natural projections.
One crucial case is when $G$ acts on itself by left translations; the resulting invariant
Brauer group is denoted by $\Br_G(G)$.

Cao used the invariant Brauer group to study local-to-global principles, weak and strong approximation
on varieties with an action of a connected linear group, see \cite{C18, C, YC}.
For a connected group $G$ he proved that 
$\Br_G(G)=\Br_1(G)$ if $G$ is {\em linear} and $\mathrm{char}(k)=0$, see \cite[Lemme 3.6]{C18}.
This naturally leads to a question raised by Cao in \cite[Exemple 4.1 (3)]{C}: 
what happens when $G=A$ is an abelian variety?  Namely, does the inclusion $\Br_A(A)\subset \Br_1(A)$ or the reverse inclusion $\Br_1(A)\subset \Br_A(A)$ hold for all abelian varieties?

In this paper we give a negative answer to the question of Cao
by showing that neither inclusion holds in general. We aim to shed light
on the properties of the somewhat mysterious invariant Brauer group
of an abelian variety. The elusive nature of $\Br_A(A)$ becomes apparent already
over an algebraically closed field $k$ of characteristic zero.
Let $A^\vee$ be the dual abelian variety of $A$, let $n\geq 2$
and let $e_{n,A}(x,y)$ be the Weil pairing $A[n]\times A^\vee[n]\to\mu_n$.
It is known \cite[p.~492]{SZ08} that if $n$ is odd, then every morphism of abelian varieties
$u\colon A\to A^\vee$ that induces a symmetric map on $n$-torsion points
$A[n]\to A^\vee[n]$, or, equivalently, such that $e_{n,A}(x,ux)=0$ for any $x\in A[n]$,
is congruent modulo $n$ to a symmetric morphism $A\to A^\vee$.
This fails if $n$ is even. In fact, $\Br_A(A)$ precisely 
measures the failure of an element $\alpha\in\Hom(A,A^\vee)/2$
such that $e_{2,A}(x,\alpha x)=0$ for any $x\in A[2]$
to come from a symmetric morphism $A\to A^\vee$. See Section \ref{tri},
in particular, formula (\ref{mrak2}) and Remark \ref{00}.

Let $\Br_A(A)(p')\subset\Br_A(A)$ be the subgroup consisting of the elements
not divisible by $\mathrm{char}(k)$.
When $k$ is algebraically closed of characteristic not equal to 2
we prove that  $\Br_A(A)(p')$ is a finite elementary 2-group, which is invariant under isogenies of odd degrees.
We give an upper bound for the order of $\Br_A(A)(p')$ in terms of the rank of the endomorphism algebra of $A$ and the Picard number of $A$. Using different techniques we describe classes of abelian varieties 
with trivial invariant Brauer group.

\bthe \label{a}
Let $A$ be an abelian variety over an algebraically closed field of characteristic different from $2$.
Then $\Br_A(A)(p')\cong(\Z/2)^n$, where 
$$n \le \rk_\Z(\End(A))-\rk_\Z(\NS(A)).$$
Moreover, $\Br_A(A)(p')=0$ in the following cases:

{\rm (i)} $A$ is a product of elliptic curves;

{\rm (ii)} $\mathrm{char}(k)=0$ and $A$ is isogenous to a power of a CM elliptic curve;

{\rm (iii)} $\mathrm{char}(k)=0$ and
$A$ is a simple abelian variety of CM type with complex multiplication by the ring of integers
of a cyclotomic field;

{\rm (iv)} $\mathrm{char}(k)>2$ and $A$ is a supersingular abelian variety.
\ethe

The upper bound is obtained in Proposition \ref{2.12}, part (i) is proved in Corollary~\ref{c1},
part (ii) is Corollary \ref{c2}, part (iii) is Proposition \ref{cyclo}, and part (iv) is Example~\ref{ex1}. 
If $E$ is an elliptic curve
over any algebraically closed field, then $\Br(E)=0$, so $\Br_A(A)\not=0$ implies $\dim(A)\geq 2$.

The question of the non-triviality of $\Br_A(A)$ over an algebraically closed field
turns out to be rather subtle.
We address this question also for complex tori, for which
the definition of the invariant Brauer group also makes sense.

\bthe \label{b}
{\rm (i)} There exist simple, non-algebraisable complex tori
in every dimension $g\geq 3$ with invariant Brauer group $\Z/2$.

{\rm (ii)} There exist simple abelian varieties over $\C$ in every dimension $g\geq 3$ 
with invariant Brauer group $\Z/2$.

{\rm (iii)} There exist abelian surfaces of CM type over $\C$, both simple and non-simple,
with invariant Brauer group $\Z/2$.
\ethe

Part (i) is proved in Section \ref{3.2}, part (ii) is proved in Section \ref{3.3}, and part
(iii) is proved in Sections \ref{surfaces} and \ref{4.2}. 
The non-simple abelian surface $A$ with $\Br_A(A)\cong\Z/2$
that we construct in Section \ref{surfaces} carries a principal polarisation, see Remark \ref{ppas}.

Assume that the ground field $k$ is an algebraically closed field 
of characteristic different from $2$. Over $k$, the invariant Brauer group of a product of abelian varieties
$A\times B$ is naturally isomorphic to $\Br_A(A)\oplus\Br_B(B)$ (up to $p$-torsion in
characteristic $p$, see Corollary \ref{c1}). We
deduce that for any integer $n\geq 0$ there is a complex abelian variety $A$ such that
$\Br_A(A)\cong(\Z/2)^n$. It would be interesting to determine the smallest dimension $g=g(n)$
for which there exists a complex abelian variety $A$ with invariant Brauer group $(\Z/2)^n$.
For example, we prove that if $g=2$, then $n\leq 2$, but 
we do not know if there is an
abelian surface over $\C$ (or another algebraically closed field) with invariant Brauer group
$(\Z/2)^2$. 

The invariant Brauer group of an abelian variety over $k$ 
is preserved by isogenies of odd degree (up to $p$-torsion in
characteristic $p$, see Proposition \ref{isogOdd}). This implies that
principally polarised $g$-dimensional abelian varieties $A$ over $\C$ such that $\Br_A(A)\not=0$ 
are dense in the complex topology of the moduli space
$\mathcal{A}_g(\C)$ for every $g\geq 2$, see Remark \ref{2.11}.
However, for a complex abelian variety $A$ which is general in the sense that
$\End(A)\cong\Z$ we have $\Br_A(A)=0$. The same holds for a general complex
torus, that is, a complex torus $A$ such that $\End(A)\cong\Z$ and $\Hom(A,A^\vee)=0$.

Finally, in Section \ref{non-closed} we study the invariant Brauer group of an abelian variety $A$ over 
a non-closed field $k$.
Let $\Br_{a,A}(A)$ be the intersection of the following three subgroups of $\Br(A)$: the invariant
subgroup $\Br_A(A)$, the algebraic subgroup $\Br_1(A)$ and the kernel of the evaluation map at the 
unit element of $A$. 

\bthe 
Let $A$ be an abelian variety over a field $k$. Then
$\Br_{a,A}(A)$ is canonically isomorphic
to the kernel of the composition of the natural maps
$$\H^1(k,\Pic(A_\s)) \lra \H^1(k,\NS(A_\s))\lra\H^1(k,\Hom(A_\s,A_\s^\vee)).$$
\ethe
This is proved in Proposition \ref{rain}. 
Using this it is easy to construct an abelian surface for which $\Br_1(A)$ is not contained in $\Br_A(A)$,
see Example \ref{end}.

\medskip

The fourth named author (Y.Z.) was partially supported by Simons Foundation
Collaboration grant $\# 585711$. Part of this work was done during his stay in
December 2019 and January 2020
at the Weizmann Institute of Science, Department of Mathematics,
whose hospitality and support are gratefully acknowledged.

The authors are grateful to the referee for the careful reading of the paper and helpful
suggestions.

\section{Definitions} \label{defi}

Let $k$ be a base field,
let $X$ be a $k$-scheme, let $G$ be a
group $k$-scheme, and let $m:G\times X\to X$ be a $k$-action of $G$ on $X$.
Write $\pi_1:G\times X\to G$ and $\pi_2:G\times X\to X$ for the natural projections.
Let $e:\Spec(k)\to G$ be the identity element of $G$. 
It induces a homomorphism $e^*:\Br(G)\to\Br(k)$, and we define
$$\Br_e(G)\colon=\Ker[e^*:\Br(G)\to\Br(k)].$$
Recall that $\Br_0(X)$ denotes the image of the natural map $\Br(k)\to\Br(X)$; we have
$\Br(G)=\Br_0(G)\oplus\Br_e(G)$.
There are induced homomorphisms
$$m^*\colon \Br(X)\to\Br(G\times X), \quad
\pi_1^*\colon \Br(G)\to\Br(G\times X), \quad
\pi_2^*\colon  \Br(X)\to\Br(G\times X).$$
The original Yang Cao's definition of the invariant Brauer group \cite[D\'ef.~3.1]{C18},
\cite[D\'ef.~1.1 (2)]{YC} 
 is as follows:
$$\Br_G(X)\colon=\{x\in\Br(X)|m^*(x)-\pi_2^*(x)\in\pi_1^*\Br(G)\}.$$
It is clear that $\Br_0(X)$ is contained in $\Br_G(X)$. 

The definition of $\Br_G(X)$ can be put in a more symmetric form.

\ble \label{2.1a}
We have $\Br_G(X)=\{x\in\Br(X)|m^*(x)\in\pi_1^*\Br(G)+\pi_2^*\Br(X)\}.$
\ele
{\em Proof.} One inclusion is trivial. To prove the other inclusion
take any $x$ in $\Br(X)$ such that $m^*(x)=\pi_1^*(y)+\pi_2^*(z)$, where $y\in\Br(G)$
and $z\in\Br(X)$. 
Modifying $y$ and $z$ by an element from the image of $\Br(k)$
we can assume without loss of generality that $y\in\Br_e(G)$.
Precomposing $m$ with $i=(e,\id)\colon X\to G\times X$ 
gives the identity map on $X$.
Hence $x=i^*(m^*(x))=e^*(y)+z=z$, so that $m^*(x)-\pi_2^*(x)\in \pi_1^*\Br(G)$. $\Box$

\medskip

In this paper we are interested in the case when $X=G$ and $m:G\times G\to G$ is the group law of $G$.
It is clear that $\Br_G(G)$ is functorial in $G$: a morphism of group $k$-schemes $G\to H$ gives rise
to a homomorphism of abelian groups $\Br_H(H)\to\Br_G(G)$.

Define
$$\Br_{e,G}(G)\colon=\Br_e(G)\cap\Br_G(G).$$
If $k$ is such that $\Br(k)=0$, e.g. $k$ is separably closed, then $\Br_G(G)=\Br_{e,G}(G)$.

\bpr \label{nuit}
Let $G$ be a group $k$-scheme with the group law $m:G\times G\to G$.
Then we have $\Br_G(G)=\Br_0(G)\oplus\Br_{e,G}(G)$ and $\Br_{e,G}(G)=\Br_e(G)\cap\Ker(\delta)$, where $$\delta=m^*-\pi_1^*-\pi_2^* : \quad\Br(G)\lra\Br(G\times G).$$
Let $[-1]\colon G\to G$ be the morphism given by inversion in $G$. Then for every
$x$ in $\Br_{e,G}(G)$ we have $[-1]^*(x)=-x$.
\epr
{\em Proof.} The first claim is clear.
Let $i_1=(\id,e)$ and $i_2=(e,\id)$ be the obvious maps $G\to G\times G$.
It is clear that on $\Br(G)$ the map $i_1^*\pi_1^*=i_1^*m^*=i_2^*\pi_2^*=i_2^*m^*$ is the identity map, whereas $i_1^*\pi_2^*=i_2^*\pi_1^*$ is the composition of $e^*\colon\Br(G)\to\Br(k)$ with
the natural injective map $\Br(k)\to\Br(G)$. By Lemma \ref{2.1a},
$x\in\Br_G(G)$ if and only if $m^*(x)=\pi_1^*(y)+\pi_2^*(z)$ for some
$y,z\in\Br(G)$. Then $x=y+e^*(z)=e^*(y)+z$, which implies
$e^*(x)=e^*(y)+e^*(z)$. When $x\in\Br_e(G)$, so that $e^*(y)+e^*(z)=0$, we can
replace $y$ by $y-e^*(y)$ and replace $z$ by $z-e^*(z)$ and hence
assume without loss of generality that $y,z\in\Br_e(G)$. Then
$x=y=z$ hence $\delta(x)=0$. Since $\Ker(\delta)$ is clearly contained in 
$\Br_G(G)$, the second claim follows.

Let $\alpha\colon G\to G\times G$ be the morphism $\alpha(x)=(x,[-1]x)$.
Then the composition 
$\alpha^*\circ\delta\colon\Br_e(G)\to\Br_e(G\times G)\to\Br_e(G)$ equals $-1-[-1]^*$,
giving the last claim. $\Box$


\medskip

We finish this section with the definition of the invariant Brauer group of a complex torus.
Recall that the cohomological Brauer group $\Br(Y)$ of a complex manifold $Y$ is defined as 
the torsion subgroup of the analytic cohomology group $\H^2(Y,\O_Y^*)$. 
If $Y$ is the complex analytic space of a complex projective variety, then $\Br(Y)$ is canonically isomorphic
to the Brauer--Grothendieck group defined using \'etale cohomology \cite[Prop.~1.3]{Sch}.

By definition, the invariant Brauer group of a complex torus $T$ is 
$$\Br_T(T):=\{x\in\Br(T)|m^*(x)\in\pi_1^*\Br(T)+\pi_2^*\Br(T)\},$$ with the same notation as 
above. In particular,
if $A$ is an abelian variety over~$\C$ and $T=A(\C)$ is the corresponding complex torus,
then $\Br_A(A)$ is canonically isomorphic to $\Br_T(T)$. The same proof as above gives the following
analogue of Proposition \ref{nuit}.

\bpr \label{torus2.2}
Let $T$ be a complex torus. Then $\Br_T(T)\subset\Br(T)$ is the kernel of 
$\delta=m^*-\pi_1^*-\pi_2^*\colon\Br(T)\to\Br(T\times T)$. 
The involution $[-1]\colon T\to T$ acts on $\Br_T(T)$ as $-1$.
\epr

\section{Abelian varieties over algebraically closed~fields} \label{tri}

\subsection{Basic properties}

\bcon \label{conv} {\rm
In this section $k$ is an algebraically closed field of $\mathrm{char}(k)\not=2$.
Let $A$ be an abelian variety over $k$. 
Let $\Br(A)(p')$ be the subgroup of $\Br(A)$ consisting of the elements of order
not divisible by $\mathrm{char}(k)$. Define
$$\Br_A(A)(p')=\Br(A)(p')\cap\Br_A(A).$$
For a commutative ring $R$ and an $R$-module $M$
we write $\wedge^2 M$ for the subgroup of $M^{\otimes 2}=M\otimes_R M$
generated by $x\otimes y-y\otimes x$ for all $x,y\in M$. We denote by $\iota$ the
inclusion $\wedge^2 M\hookrightarrow M^{\otimes 2}$. This will be applied when $R$ is
$\Z$, $\Z_2$ or $\Z/2$.}
\econ

\brem {\rm
The discussion and results of this section can be easily modified for the case of complex tori by considering 
classical cohomology groups instead of \'etale cohomology.}
\erem

The case of  dimension~1 does not present any difficulty.

\brem \label{2.1} {\rm
For an elliptic curve $E$ over any algebraically closed field we have $\Br(E)=0$, 
see \cite[III, Cor. 1.2]{Gro68}. In particular, $\Br_E(E)=0$.}
\erem

Using the Kummer exact sequence it is easy to show that
the  multiplication by $n$ map $[n]\colon A\to A$
induces multiplication by $n^2$ on $\Br(A)(p')$, see \cite[Sect.~2]{Ber72}. 
In particular, $[-1]\colon A\to A$ acts on $\Br(A)(p')$ trivially.
This fact and the last statement of Proposition \ref{nuit} imply
that $\Br_A(A)(p')$ is contained in $\Br_A(A)[2]$. By Proposition \ref{nuit} we conclude that $\Br_A(A)(p')$
is the kernel of 
$$
\delta\colon \ \Br(A)[2]\lra \Br(A\times A)[2].
$$
The long exact sequence of \'etale cohomology attached to the
Kummer exact sequence allows us to write this map as
$$
\delta\colon \ \frac{\H_\et^2(A,\Z/2)}{\NS(A)/2} \lra \frac{\H_\et^2(A\times A,\Z/2)}{\NS(A\times A)/2}.
$$
To compute this map we need to recall some standard facts about abelian varieties over an algebraically
closed field.

Let $\si\colon A\times A\to A\times A$ be the involution that exchanges the factors. By abuse of notation we 
also denote by $\si$ the maps induced by $\si$ on cohomology groups of $A\times A$. 
Recall that $m\colon A\times A\to A$
is the multiplication map, and $\pi_1\colon A\times A\to A$ and $\pi_2\colon A\times A\to A$ are the 
natural projections. As usual, $A^\vee$ is the dual abelian variety of $A$. 
The dual morphism $m^\vee\colon A^\vee\to A^\vee\times A^\vee$ is the diagonal map.

Let $\H_\et^i(-,\Z_2)$ be the 2-adic \'etale cohomology groups.
Since $A$ is an abelian variety, the cup-product 
$$\cup_A \colon \H_\et^1(A,\Z_2)\times \H_\et^1(A,\Z_2)\lra \H_\et^2(A,\Z_2)$$
gives rise to a canonical isomorphism of $\Z_2$-modules
$$\wedge^2 \H_\et^1(A,\Z_2)\cong \H_\et^2(A,\Z_2).$$
The K\"unneth decomposition in degree 1 is
\begin{equation}
\H_\et^1(A\times A,\Z_2)\cong \H_\et^1(A,\Z_2)\oplus \H_\et^1(A,\Z_2), \label{kun1}
\end{equation}
where the direct summands are embedded via the induced maps $\pi_1^*$ and $\pi_2^*$.
Precomposing $m$ with $A\stackrel{\sim}\lra A\times \{0\}\to A\times A$ gives the identity map,
hence it is clear that $m^*\colon\H_\et^1(A,\Z_2)\to \H_\et^1(A\times A,\Z_2)$ sends $x$ to $x\oplus x$.
Taking the second exterior power in (\ref{kun1}) we obtain
the K\"unneth decomposition in degree 2:
\begin{equation}
\H^2_\et(A\times A,\Z_2)\cong
\H_\et^2(A,\Z_2)\oplus\H_\et^2(A,\Z_2)\oplus \left(\H_\et^1(A,\Z_2)\otimes_{\Z_2} \H_\et^1(A,\Z_2)\right),
\label{ku}
\end{equation}
where the first factor is embedded via $\pi_1^*$,
the second factor via $\pi_2^*$, and the third factor 
via the map that sends $x\otimes y$ to
$ \pi_1^*(x)\cup_{A\times A}\pi_2^*(y)$. 
The cup-product
$$\cup_{A\times A}\colon \ \ \H_\et^1(A\times A,\Z_2)\otimes_{\Z_2} \H_\et^1(A\times A,\Z_2)\lra 
\H_\et^2(A\times A,\Z_2)$$
is skew-symmetric, hence $\si^*$ sends $\pi_1^*(x)\cup_{A\times A}\pi_2^*(y)$ to
$$\pi_2^*(x)\cup_{A\times A}\pi_1^*(y)=-\pi_1^*(y)\cup_{A\times A} \pi_2^*(x).$$ It follows that
$\si$ acts on the third summand of (\ref{ku}) by sending $x\otimes y$ to $-y\otimes x$.
 
Taking the cup product of $m^*\colon \H_\et^1(A,\Z_2) \to \H_\et^1(A\times A,\Z_2)$ with itself, we find that
in terms of (\ref{ku}) the map $m^*\colon\H_\et^2(A,\Z_2)\to \H_\et^2(A\times A,\Z_2)$ is 
$(\id,\id,\iota)$, where $\iota$ is defined in Convention \ref{conv}. 
Thus $\delta\colon \H_\et^2(A,\Z_2) \to \H_\et^2(A\times A,\Z_2)$ is the map $(0,0,\iota)$. 

Write $V=\H_\et^1(A,\Z_2)$. This is a free $\Z_2$-module of rank $2g$, where $g=\dim(A)$.
Using the Kummer sequence and the Weil pairing we obtain canonical isomorphisms
\begin{equation}
V/2\cong\H_\et^1(A,\Z/2)\cong A^\vee[2]\cong\Hom(A[2],\Z/2). \label{V2}
\end{equation}
We have $\wedge^2 V=(V^{\otimes 2})^{\si}\cong \H_\et^2(A,\Z_2)$.

Applying the canonical isomorphism $\NS(A)\cong\Hom(A,A^\vee)^{\rm sym}$ 
to $A\times A$ instead of $A$, we obtain
\begin{equation}
\NS(A\times A)\cong\NS(A)\oplus\NS(A)\oplus\Hom(A,A^\vee); \label{ns}
\end{equation}
see also \cite[Prop.~1.7]{SZ14} or Lemma \ref{known} below. The map $m^*\colon\NS(A)\to\NS(A\times A)$
is obtained from the map $\Hom(A,A^\vee)\to \Hom(A\times A, A^\vee\times A^\vee)$
that sends $\phi$ to $m^\vee\circ\phi\circ m$. Using that $m^\vee$ is the diagonal map,
we deduce that in terms of (\ref{ns}) the map $\delta=m^*-\pi_1^*-\pi_2^*$ is 
the natural inclusion $\Hom(A,A^\vee)^{\rm sym}\hookrightarrow \Hom(A,A^\vee)$.
The involution
$\si$ acts on $\NS(A\times A)$ by permuting the first two direct summands
in (\ref{ns}) and preserving the third one, on which it acts by $\phi\mapsto\phi^\vee$. 
The 2-adic first Chern class map sends
each direct summand in (\ref{ns}) to
the corresponding direct summand in (\ref{ku}). 
For the last summand the mod 2 first Chern class gives a map
\begin{equation}
\Hom(A,A^\vee)/2\hookrightarrow (V/2)^{\otimes 2}\cong\Hom(A[2],A^\vee[2]) \label{hom2}
\end{equation}
where the isomorphism $(V/2)^{\otimes 2}\cong\Hom(A[2],\Z/2)\otimes_\Z A^\vee[2]\cong\Hom(A[2],A^\vee[2])$ comes from (\ref{V2}).

\brem \label{0}{\rm At least
when ${\rm char}(k)=0$, this map is induced by the natural action of homomorphisms of 
abelian varieties on torsion points.
(If $2$ is replaced by an arbitrary positive integer, the analogous map is the
negative of the action on torsion points.) This is proved by applying \cite[Lemma 2.6]{OSZ} 
to $A\times A$.}
\erem

Putting together the descriptions of $\delta$ for the cohomology groups and the N\'eron--Severi groups
we obtain that $\Br_A(A)(p')$ is the kernel of the natural map
\begin{equation}
\frac{\wedge^2 (V/2)}{\NS(A)/2} \lra \frac{(V/2)^{\otimes 2}}{\Hom(A,A^\vee)/2}.
\label{mrak}
\end{equation}
In order to understand the definition of this map, we recall that $\wedge^2 (V/2)$ and $\Hom(A,A^\vee)/2$ are subgroups of $(V/2)^{\otimes 2}$, using the map (\ref{hom2}) for $\Hom(A,A^\vee)/2$.

We remind the reader that $\wedge^2 V=(V^{\otimes 2})^{\si}$ and 
$\NS(A)=\Hom(A,A^\vee)^\si$.
Write $$L=\Hom(A,A^\vee)\otimes\Z_2.$$
This is a primitive $\Z_2$-sublattice of $V^{\otimes 2}$ by the proof of \cite[Lemma~1]{Tat66}. 
The involution $\si$ acts on $L$ and we have
$$\NS(A)\otimes\Z_2=L^{\si}=L\cap \wedge^2 V\subset V^{\otimes 2}.$$ 
In the following proposition, we write $(L/2)\cap\, \wedge^2 (V/2)$ 
(respectively, $(L/4)\cap\, \wedge^2 (V/4)$) for the intersection in $(V/2)^{\otimes 2}$ 
(respectively, in $(V/4)^{\otimes 2}$).

\bpr \label{uno}
We have 
$\Br_A(A)(p')=\frac{(L/2)\ \cap \ \wedge^2 (V/2)}{(L \, \cap \, \wedge^2 V)/2}
=\frac{(L/2)\ \cap \ \wedge^2 (V/2)}{((L/4) \, \cap \, \wedge^2 (V/4))/2}.$
\epr
{\em Proof.}
The first equality is clear, so let us prove that the image of 
$L \, \cap \, \wedge^2 V$ in $(L/2)\ \cap \ \wedge^2 (V/2)$ is equal to the (a priori, larger) image of 
$(L/4) \, \cap \, \wedge^2 (V/4)$. Take any $x\in (L/4)\cap \wedge^2 (V/4)$ and lift it 
to $\tilde x\in L$. Clearly, $\tilde x+\si(\tilde x)\in L^\si=L \cap \wedge^2 V$.
Since $\tilde x$ and $\si(\tilde x)$ are congruent modulo 2, we can write 
$\tilde x+\si(\tilde x)=2v$ for some $v\in L^\si$. Furthermore, $\tilde x$ and $\si(\tilde x)$ 
are congruent modulo 4, so $v$ is congruent to $\tilde x$ modulo 2, that is,
$v\in L \cap \wedge^2 V$ is a lifting of the image of $x$ in $(L/2)\ \cap \ \wedge^2 (V/2)$. $\Box$

\bpr \label{duo}
Let $\Z/2$ be the group generated by $\si$. Then $\Br_A(A)(p')$ is 
the kernel~of 
\begin{equation}
\label{LtoV2}
\H^1(\Z/2,L) \lra \H^1(\Z/2, V^{\otimes 2})=V/2\cong (\Z/2)^{2g},
\end{equation}
where $\si$ acts on $L=\Hom(A,A^\vee)\otimes\Z_2$ by $\phi\mapsto\phi^\vee$
and on $V^{\otimes 2}$ by $x\otimes y\mapsto-y\otimes x$.
\epr
{\em Proof.} Using (\ref{mrak}) we identify $\Br_A(A)(p')$ with the kernel of the natural map
\begin{equation}
\frac{\wedge^2 (V/2)}{L^\si/2} \lra \frac{((V/2)^{\otimes 2})^\si}{(L/2)^\si}.
\label{mrak1}
\end{equation}
Alternatively, $\Br_A(A)(p')$ is the kernel of the natural map
\begin{equation}
\frac{(L/2)^\si}{L^\si/2} \lra \frac{((V/2)^{\otimes 2})^\si}{\wedge^2 (V/2)}=
\frac{S^2 (V/2)}{\wedge^2 (V/2)}\tilde\longleftarrow V/2,
\label{mrak2}
\end{equation}
where the last map is given by $x\mapsto x\otimes x$.

For any abelian group $M$ with trivial 2-torsion, we have a short exact sequence
$$ 0 \lra M \overset{[2]}{\lra} M \lra M/2 \lra 0. $$
If $M$ is acted on by $\Z/2$ with generator $\si$, then the associated long exact sequence in cohomology leads to a further exact sequence
$$0\lra M^\si/2\lra (M/2)^\si\lra \H^1(\Z/2,M)\lra 0.$$
This gives the interpretation in terms of $\H^1(\Z/2,-)$ in (\ref{LtoV2}).
$\Box$

\brem \label{00} {\rm This remark uses Remark \ref{0} so we assume ${\rm char}(k)=0$.
Via the natural identification $V/2 \cong \Hom(A[2], \Z/2)$,
the map in (\ref{mrak2}) sends
$u\in (\Hom(A,A^\vee)/2)^\si$ to the linear form on $A[2]$ whose value on
$x\in A[2]$ is $e_{2,A}(x,ux)$, where $e_{2,A}(x,y)$ is the Weil pairing 
$A[2]\times A^\vee[2]\to\Z/2$.}
\erem

\bex \label{ex1} {\rm
Suppose that $\mathrm{char}(k)>2$ and $A$ is a {\em supersingular} abelian variety \cite{LO98}.
Then it is well known that $\End(A)$ is a free $\Z$-module of rank $(2g)^2$. Since $A$ and $A^{\vee}$ are isogenous,
$\Hom(A,A^{\vee})$ is also a free $\Z$-module of rank $(2g)^2$. This implies that
$L$ is a free $\Z_2$-module of rank $(2g)^2$. Taking into account that $L$ is a primitive submodule in $V^{\otimes 2}$,
and the latter is also a  free $\Z_2$-module of the same rank $(2g)^2$, we conclude that $L=V^{\otimes 2}$. 
This implies  that  \eqref{LtoV2} is the identity map, so by Proposition \ref{duo} we have
$\Br_A(A)(p')=0$ if $A$ is supersingular.}
\eex

\bco \label{c1}
For abelian varieties $A$ and $B$ there is a canonical isomorphism
$$\Br_{A\times B}(A\times B)(p')\cong\Br_A(A)(p')\oplus\Br_B(B)(p').$$
In particular, if $A$ is the product of elliptic curves, then $\Br_A(A)(p')=0$.
\eco
{\em Proof.} We have
$$\Hom(A\times B,A^\vee\times B^\vee)=\big(\Hom(A,A^\vee)\oplus \Hom(B,B^\vee)\big)\oplus
\big(\Hom(A,B^\vee)\oplus \Hom(B,A^\vee)\big).$$
Here the second summand on the right hand side is a permutation $\Z/2$-module, hence it has
trivial $\H^1(\Z/2,-)$ group. The displayed formula now follows from Proposition \ref{duo}.
The second statement follows from Remark \ref{2.1}. $\Box$

\bpr \label{isogOdd}
An isogeny of abelian varieties of odd degree $\phi: A \to B$ induces an isomorphism
$\Br_A(A)(p') \cong \Br_B(B)(p')$.
\epr
{\em Proof.}
Let $m=\deg(\phi)$.
There is an isogeny $\psi\colon B \to A$ such that
$$\psi\circ \phi=[m]_A, \quad\quad  \phi\circ \psi=[m]_B.$$
By functoriality,  $\phi$ and $\psi$ induce homomorphisms 
$$\phi^*\colon \Br(B) \lra \Br(A), \quad\quad \psi^*\colon \Br(A) \lra \Br(B).$$
On $\Br(A)(p')$, the composition $\phi^*\psi^*$ is the multiplication by $m^2$.
Similarly, on $\Br(B)(p')$, the composition  $\psi^*\phi^*$ is the multiplication by $m^2$.
We have
$$\phi^*(\Br_B(B)(p'))\subset \Br_A(A)(p'), \quad\quad \psi^*(\Br_A(A)(p'))\subset \Br_B(B)(p').$$
Since $m$ is odd and both $\Br_A(A)(p')$ and $\Br_B(B)(p')$ are annihilated by 2, the map $\phi^*$
induces an isomorphism $\Br_B(B)(p')\stackrel{\sim}\lra\Br_A(A)(p')$. $\Box$

\brem  {\rm
Recall \cite[pp. 14--15]{LO98} that if $\mathrm{char}(k)=p>2$ and  $A$ is a {\em supersingular} abelian variety then there is a {\em minimal isogeny}
$\phi\colon A \to B$ such that $\deg(\phi)$ is a power of $p$ and $B$ is a product of (supersingular) elliptic curves. Now Proposition \ref{isogOdd}
combined with Corollary \ref{c1} gives another proof that $\Br_A(A)(p')=0$.}
\erem

\bco \label{c2}
Suppose that $\mathrm{char}(k)=0$.
Let $A$ be an abelian variety isogenous to a power of an elliptic curve
with complex multiplication. Then $\Br_A(A)=0$.
\eco
{\em Proof.} By a result of C. Schoen \cite{Sch92} in this case
$A$ is isomorphic to a product of simple abelian varieties,
which necessarily are elliptic curves. 
We conclude by applying Corollary \ref{c1}. $\Box$

\brem \label{2.10} {\rm 
If $\End(A)\otimes\Q$ is a direct sum of totally real fields,
for example, when $A$ is simple and
has Albert type I, then $\si$ acts trivially on $\End(A)$
hence also on $\Hom(A,A^\vee)$, so $\Br_A(A)(p')=0$. In particular, for a general abelian
variety (when $\End(A)=\Z$) we have $\Br_A(A)(p')=0$. Similarly, for a general complex
torus (when $\Hom(A,A^\vee)=0$) we have $\Br_A(A)=0$.}
\erem

\brem \label{2.11} {\rm
Let $\mathcal{A}_g$ be the coarse moduli space of 
principally polarised abelian varieties of dimension $g \ge 2$ over $\C$. 
The symplectic group $\Sp(2g,\R)$ acts transitively 
on the Siegel upper half-plane $\mathfrak{H}_g$, and
we have $\mathcal{A}_g(\C)=\Sp(2g,\Z)\setminus \mathfrak{H}_g$.

 Let $(A,\lambda)$ and $(B,\mu)$ be complex principally polarised $g$-dimensional abelian varieties. One says that $(A,\lambda)$ and $(B,\mu)$ are in the same {\em Hecke orbit} if there exist a $g$-dimensional   complex abelian variety $C$ and isogenies
$\alpha: C \to A$ and $\beta: C \to B$
such that the polarisations
$\alpha^{*}\lambda$ and $\beta^{*}\mu$ of $C$ are equal \cite[pp.~441--442]{CO09}.
If neither $\deg(\alpha)$ nor
$\deg(\beta)$ is divisible by a prime $p$, then one says that  $(A,\lambda)$ and $(B,\mu)$ are
 in the same {\em prime-to-$p$  Hecke orbit}. It follows from Proposition \ref{isogOdd} that if
 $(A,\lambda)$ and $(B,\mu)$ are in the same prime-to-$2$ Hecke orbit, then
 $\Br_A(A)\cong \Br_B(B)$.

 It is known \cite[p.~92]{Ch05} that each prime-to-$p$ Hecke orbit is dense in the complex topology of  $\mathcal{A}_g(\C)$. Indeed, consider the group scheme $G=\Sp(2g)$ over $\Z$.
The connected semisimple algebraic group $G_\Q$ over $\Q$ is
simply connected so it satisfies weak approximation \cite[Prop.~7.9]{PR}. In particular, $G(\Q)$
is dense in $G(\R)\times G(\Q_p)$, and since $G(\Z_p)$ is open in $G(\Q_p)$, we see that
$G(\Q)\cap G(\Z_p)$ is dense in $G(\R)$. But a prime-to-$p$  Hecke orbit in $\mathcal{A}_g$
is the image of a $G(\Q)\cap G(\Z_p)$-orbit in $\mathfrak{H}_g$. 
(Moreover, a prime-to-$p$  Hecke orbit is equidistributed in $\mathcal{A}_g$
with respect to the normalised Haar measure, as follows from
equidistribution of Hecke points in $\Sp(2g,\Z)\setminus\Sp(2g,\R)$, see \cite[Thm.~1.6]{COU01}.)
 
  Taking $p=2$ and applying Corollary \ref{c1} together with Remark \ref{ppas}
below, we see that  the set of (isomorphism classes of) principally polarised complex $g$-dimensional abelian varieties $A$ such that $\Br_A(A)\not=0$ is dense in the
complex topology of $\mathcal{A}_g(\C)$.  The same applies to abelian varieties $A$ with $\Br_A(A)=0$.}
\erem


Proposition \ref{duo} has an analogue for complex tori, with the same proof.

\bpr \label{torus3.5}
Let $T$ be a complex torus. Then $\Br_T(T)$
is the kernel of the map analogous to the map $(\ref{LtoV2})$, where $L=\Hom(T,T^{\vee})\otimes\Z_2$ and 
$V=\H^1(T,\Z_2)=\H^1(T,\Z)\otimes_\Z\Z_2$.
\epr

\subsection{Upper bound} \label{bound}

Write $H = \Hom(A,A^\vee)\otimes\Q$ and $H^- = \{ \phi \in H | \phi^\vee=-\phi \}$.
Note that the choice of a polarisation of $A$ induces an isomorphism of $\Q$-vector spaces $H \cong \End(A)\otimes\Q$ under which $\phi\mapsto\phi^\vee$ corresponds to the Rosati involution,
see \cite[Ch.~20]{M}.

The following proposition
gives an upper bound for the size of $\Br_A(A)(p')$.  
The right hand side of this inequality is isogeny-invariant.

\bpr \label{2.12}
$\dim_{\F_2}(\Br_A(A)(p')) \leq \dim_\Q(H^-)=\rk_\Z(\End(A))-\rk_\Z(\NS(A))$.
\epr
{\em Proof.} From (\ref{mrak2}), we have
$$ \Br_A(A)(p') \subset \frac{(L/2)^\si}{L^\si/2} \subset \frac{L/2}{L^\si/2}. $$
Since $L^\si$ is primitive in $L$, we have
$$ \dim_{\F_2} \left( \frac{L/2}{L^\si/2} \right) = \rk_{\Z_2}(L) - \rk_{\Z_2}(L^\si) = \dim_\Q(H) - \dim_\Q(H^\si) = \dim_\Q(H^-). $$
It is clear that $\dim_\Q(H^-)=\rk_\Z(\End(A))-\rk_\Z(\NS(A))$. $\Box$

\medskip

The following lemma is needed for the proof of Proposition \ref{maxdim} below.
It is due to G.~Shimura
\cite[Prop.~15, p.~177]{Shi63}, see also \cite[Subsection (4.1), p.~488]{O88}. 

\ble
\label{shimura}
Let $A$ be a simple positive-dimensional abelian variety over an algebraically closed field $k$ of characteristic $0$. Suppose that $A$ is of type III in Albert's classification, that is,
$\End^0(A):=\End(A)\otimes \Q$ is a totally definite quaternion algebra over a totally real number field $F$.
Then we have 
$$2\,\dim(A)=m\, \dim_{\Q}(\End^0(A))=4m [F:\Q]=4m\, \rk(\NS(A))$$ 
for some integer $m\geq 2$.
\ele
{\em Proof.}
Withous loss of generality we may assume that $k=\C$. 
Let us choose a polarisation on $A$. The associated Rosati involution is 
the standard involution on the quaternion $F$-algebra $\End^0(A)$, hence the space of
symmetric elements of $\End^0(A)$ is $F\subset \End^0(A)$, see \cite [Ch.~21]{M}. In particular, 
$\rk(\NS(A))=[F:\Q]$.

It is known \cite[Ch.~21]{M} that $\dim_{\Q}(\End^0(A))$ divides $2\dim(A)$, so $m\in\Z$, $m\geq 1$.
If $m=1$, then \cite[Prop.~15, p.~177]{Shi63} implies that $A$ is isogenous to the square of 
a certain abelian variety. In particular, $A$ is not simple, which contradicts our assumption. 
$\Box$

\bpr \label{maxdim}
Suppose that $\mathrm{char}(k)=0$.
Let the abelian variety~$A$ be isogenous to $A_1^{n_1} \times \dotsb \times A_m^{n_m}$ where the $A_i$ are simple and pairwise non-isogenous.
Then
$ \dim_\Q(H^-) \leq \max(n_i) \dim(A) $.
\epr
{\em Proof.} The endomorphism algebra of $A$ is described by
\begin{equation}
\End^0(A) \cong \prod_{i=1}^m \Mat_{n_i}(D_i)
\label{enda}
\end{equation}
where $D_i = \End^0(A_i)$ are division algebras.

Choose a polarisation on $A$ which pulls back to a product of polarisations on the $A_i$.
The associated Rosati involution preserves the decomposition (\ref{enda}) and acts on $\Mat_{n_i}(D_i)$ as the composition of matrix transpose with the entry-wise Rosati involution of $D_i$.
Hence
$$ \dim_\Q(H^-) = \sum_{i=1}^m \Bigl( \frac{n_i(n_i-1)}{2} \dim_\Q(D_i) + n_i\,\dim_\Q(D_i^-) \Bigr).$$
The division algebras which occur as the endomorphism algebras of simple abelian varieties were classified by Albert, see \cite[Ch.~21]{M}.
From this classification, we obtain the following bounds for the dimension.

\smallskip

\setlength\extrarowheight{2pt}
\begin{tabular}{ccc}
Endomorphism type & Maximum value of $\dim_\Q(D_i)$ & $\dim_\Q(D_i^-)$
\\ I   & $\dim(A_i)$  & 0
\\ II  & $2\dim(A_i)$ & $\frac{1}{4}\dim_\Q(D_i)$
\\ III & $\dim(A_i)$  & $\frac{3}{4}\dim_\Q(D_i)$
\\ IV  & $2\dim(A_i)$ & $\frac{1}{2}\dim_\Q(D_i)$
\end{tabular}

\medskip

\noindent Most of this table is based on \cite[p.~202]{M}, noting that Mumford's $\eta$ is equal to $1-\dim_\Q(D_i^-)/\dim_\Q(D_i)$.
The only entry in this table which is not taken from \cite[p.~202]{M} is the maximum value of $\dim(D_i)$ for type III. In this case the desired result follows from Lemma \ref{shimura},
which gives us a better bound than \cite{M} for type III.

From the table, we deduce that
$$ \dim_\Q(D_i) \leq 2\dim(A_i), \quad \dim_\Q(D_i^-) \leq \dim(A_i). $$
Hence
$$ \frac{n_i(n_i-1)}{2} \dim_\Q(D_i) + n_i\,\dim_\Q(D_i^-) \leq n_i^2\dim(A_i) = n_i \dim(A_i^{n_i}). $$
Summing this over the isotypic factors of $A$ proves the proposition.
$\Box$

\bco
Suppose that $\mathrm{char}(k)=0$.
If $A$ is an abelian variety over $k$ of dimension~$2$, then $\dim_{\F_2}(\Br_A(A)) \leq 2$.
\eco
{\em Proof.} By Proposition~\ref{maxdim}, this is true whenever $A$ is simple or is isogenous to a product of non-isogenous elliptic curves.
It remains to check the case where $A$ is isogenous to the square of an elliptic curve $E$.

If $\End(E)=\Z$, then $\End^0(A) \cong \Mat_2(\Q)$ and the Rosati involution is transposition.
So $\dim_\Q(H^-) = 1$.

If $E$ has complex multiplication, then $\dim_\Q(H^-) = 4$ so we cannot deduce the corollary from Proposition~\ref{maxdim}.
In this case $\Br_A(A)=0$ by Corollary~\ref{c2}.
$\Box$

\subsection{Calculating \texorpdfstring{$\Br_A(A)$}{BrA(A)} for complex tori}

A \textit{complex structure} on a real vector space $\W$ is an element $J \in \End(\W)$ satisfying $J^2=-1$.
Giving a complex structure on $\W$ is equivalent to giving a complex vector space whose underlying real vector space is equal to $\W$, by letting $J$ represent the action of $i$ on the complex vector space.

Let $A$ be a complex torus of dimension~$g$.
By definition, $A = \V/\Lambda$ where $\V$ is a complex vector space of dimension~$g$ and $\Lambda$ is a discrete free abelian subgroup of $\V$ of rank~$2g$.
We can naturally identify $\Lambda$ with $\H_1(A,\Z)$.
We also identify $\V$ with $\Lambda_\R$ equipped with a complex structure $J$.

Let $A'$ be another complex torus with $\Lambda' = \H_1(A',\Z)$ and with associated complex structure $J' \in \End(\Lambda'_\R)$.
Then each homomorphism $A \to A'$
 induces a homomorphism of homology groups $\Lambda \to \Lambda'$.
Conversely, a homomorphism $f \colon \Lambda \to \Lambda'$ comes from a homomorphism of complex tori if and only if it intertwines the complex structures on $\Lambda_\R$ and $\Lambda'_\R$.
In other words, there is a natural bijection
\begin{equation} \label{hom-tori}
\Hom(A, A') = \{ f \in \Hom(\Lambda, \Lambda') | f \circ J = J' \circ f \}.
\end{equation}

The \textit{dual complex torus} $A^\vee$ of $A$ is the quotient of the vector space
of antilinear homomorphisms $f\colon \V\to\C$ by the subgroup 
$\Lambda^\vee=\{f|\Im f(\Lambda)\subset\Z\}$, see \cite[Ch.~9]{M}. Thus
 $\Lambda^\vee \cong \Hom(\Lambda, \Z)$ and
$A^\vee = \Lambda^\vee_\R/\Lambda^\vee$, where 
the complex structure of $\Lambda^\vee_\R$ is given by $-J^*$, see \cite[p.~5]{Kem91}.
This gives a canonical isomorphism $(\Lambda^\vee)^\vee\cong\Lambda$. Using this isomorphism,
the canonical pairing $\Lambda^\vee\times \Lambda\to\Z$
is identified with the {\em negative} of the canonical pairing $\Lambda\times\Lambda^\vee\to\Z$.

Write $\Bi(\Lambda)$ for the set of bilinear forms $\Lambda \times \Lambda \to \Z$. Let
$\tau$ be the involution of $\Bi(\Lambda)$ which exchanges the arguments of a bilinear form.

Write $\Bi_J(\Lambda)$ for the set of $J$-invariant bilinear forms, that is, bilinear forms $B \in \Bi(\Lambda)$ for which the induced form on $\Lambda_\R$ satisfies $B(x,y) = B(Jx,Jy)$ for all $x,y \in \Lambda_\R$.
Since $J^2 = -1$, this is equivalent to saying
\begin{equation} \label{BJ}
\Bi_J(\Lambda) = \{ B \in \Bi(\Lambda) | B(Jx,y) = B(x,-Jy) \text{ for all } x,y \in \Lambda_\R \}.
\end{equation}
There is a natural bijection
$\Hom(\Lambda, \Lambda^\vee) \tilde\to \Bi(\Lambda)$
which sends \( f \in \Hom(\Lambda, \Lambda^\vee) \) to \( B_f(x,y) = f(x)(y) \).
By (\ref{BJ}), we have $B_f \in \Bi_J(\Lambda)$ if and only if $f \circ J = -J^* \circ f$.

Hence, using (\ref{hom-tori}) for \( A' = A^\vee \), we obtain natural bijections
\begin{equation} \label{BiJ}
\Hom(A, A^\vee) = \{ f \in \Hom(\Lambda, \Lambda^\vee) | f \circ J = -J^* \circ f \} = \Bi_J(\Lambda).
\end{equation}
Under these bijections, the involution $\sigma$ of $\Hom(A, A^\vee)$ given by $\phi \mapsto \phi^\vee$ corresponds to $-\tau$. (This follows from
the fact that the canonical pairings $\Lambda^\vee\times \Lambda\to\Z$
and $\Lambda\times\Lambda^\vee\to\Z$ differ by sign.) 
In particular, $\NS(A) = \Hom(A, A^\vee)^{\rm sym}$ corresponds to the alternating forms in $\Bi_J(\Lambda)$.

An element $\Hom(A, A^\vee)^{\rm sym}$ is a \textit{polarisation} if and only if the corresponding alternating form $E \in \Bi_J(\Lambda)$ satisfies the following condition:
\begin{equation} \label{polarisation}
\text{the symmetric bilinear form } E(Jx,y) \text{ is positive definite}.
\end{equation}
Note that $E(Jx,y)$ is the real part of the Hermitian form attached to $E$ as in \cite[p.~19, Lemma]{M}.
Recall that a complex torus $A$ is an abelian variety if and only if there exists a polarisation of $A$ \cite[p.~35, Corollary]{M}.
A complex torus is said to be \textit{non-algebraisable} if it does not possess a polarisation.

\ble \label{bi-ends}
If $\Bi_J(\Lambda)$ contains a non-degenerate bilinear form (which is always true when $A$ is an abelian variety), then $\rk_\Z (\Bi_J(\Lambda)) = \rk_\Z (\End(A))$.
\ele
{\em Proof.}
Since $\Bi_J(\Lambda)$ and $\End(A)$ are both free $\Z$-modules, it suffices to prove that $\dim_\Q(\Bi_J(\Lambda_\Q)) = \dim_\Q(\End(A) \otimes \Q)$.

Let $B$ be a non-degenerate form in $\Bi_J(\Lambda)$.
Then every bilinear form on $\Lambda_\Q$ can be written as $(x,y) \mapsto B(x, uy)$ for some $u \in \End(\Lambda_\Q)$.
The form $B(x,uy)$ is in $\Bi_J(\Lambda_\Q)$ if and only if $u$ commutes with $J$.
The lemma is proved because $\End(A) \otimes \Q = \End_J(\Lambda_\Q)$.
$\Box$

\medskip

A bilinear form $B \colon \Lambda \times \Lambda \to \Z$ is said to be \textit{even} if $B(x,x) \equiv 0 \bmod 2$ for all $ x \in \Lambda$.
Write $\Bi_J(\Lambda)^{\mathrm{sym,even}}$ for the set of symmetric, even forms in $\Bi_J(\Lambda)$.

Write $\Lambda_2 = \Lambda \otimes \Z_2$.
Then $\Bi_J(\Lambda_2) = \Bi_J(\Lambda) \otimes \Z_2$ and $\Bi_J(\Lambda_2)^{\mathrm{sym,even}} = \Bi_J(\Lambda)^{\mathrm{sym,even}} \otimes \Z_2$.

\bpr \label{calc}
$\Br_A(A) \cong \Bi_J(\Lambda)^{\mathrm{sym,even}} / (1+\tau)\Bi_J(\Lambda)$.
\epr
{\em Proof.}
Considering the action of $\Z/2$ on $L = \Hom(A, A^\vee) \otimes \Z_2$ by $\sigma$, we have
\[ \H^1(\Z/2, L) = L^{-\sigma}/(1-\sigma)L = \Bi_J(\Lambda_2)^{\mathrm{sym}} / (1+\tau)\Bi_J(\Lambda_2). \]
Note that $V = \H^1(A,\Z_2)$ can be canonically identified with $\Lambda_2^\vee$ and hence $V^{\otimes 2}$ can be identified with the set of bilinear maps $\Lambda_2 \times \Lambda_2 \to \Z_2$.
The bijection~\eqref{BiJ} identifies the inclusion $L \subset V^{\otimes 2}$ with the inclusion $\Bi_J(\Lambda_2) \to V^{\otimes 2}$.

Consequently Proposition~\ref{torus3.5} 
tells us that $\Br_A(A)$ is isomorphic to the kernel of the natural map
\[ \Bi_J(\Lambda_2)^{\mathrm{sym}} / (1+\tau)\Bi_J(\Lambda_2) \lra V/2 \]
which sends $B \in \Bi_J(\Lambda_2)^{\mathrm{sym}}$ to the function $x \mapsto B(x,x) \bmod 2$.
In other words, the kernel of this map consists precisely of the even forms.

Since
$$ 2\Bi_J(\Lambda)^{\mathrm{sym,even}} \subset 2\Bi_J(\Lambda)^{\mathrm{sym}} \subset (1+\tau)\Bi_J(\Lambda), $$
we have
$\Bi_J(\Lambda_2)^{\mathrm{sym,even}} / (1+\tau)\Bi_J(\Lambda_2)= \Bi_J(\Lambda)^{\mathrm{sym,even}} / (1+\tau)\Bi_J(\Lambda)$.
$\Box$

\medskip

In many of our examples, we will look at abelian varieties $A$ for which $K = \End(A) \otimes \Q$ is a product of CM fields satisfying $\dim_\Q(K) = \dim(A)$.
This includes simple abelian varieties of CM type.
In this case, we can identify $\Lambda_\Q$ with $K$.
Consequently, for each $\alpha \in K$, it makes sense to define a bilinear form 
$$B_\alpha \colon \Lambda_\Q \times \Lambda_\Q \to \Q, \quad
B_\alpha(x, y) = \tr_{K/\Q}(\alpha x \bar y).$$
Let $D_K(\Lambda) = \{ \alpha \in K | B_\alpha(\Lambda \times \Lambda) \subset \Z \}$.
Note that if $\Lambda$ is a complex conjugation-invariant subring of $K$, then we can simplify the calculation further as in this case,
\begin{equation} \label{DK-subring}
D_K(\Lambda) = \{ \alpha \in K | \tr_{K/\Q}(\alpha x) \in \Z \text{ for all } x \in \Lambda \}.
\end{equation}
The bilinear forms $B_\alpha$ are $J$-invariant, so thanks to Lemma~\ref{bi-ends} we have
$$\Bi_J(\Lambda) = \{ B_\alpha | \alpha \in D_K(\Lambda) \}.$$
Note also that $\tau(B_\alpha) = B_{\iota(\alpha)}$ where $\iota$ denotes complex conjugation in $K$.
Writing $D_K(\Lambda)^{\iota,\mathrm{even}} = 
\{ \alpha \in D_K(\Lambda) | B_\alpha \text{ is symmetric and even} \}$, we conclude from Proposition~\ref{calc} that
\begin{equation} \label{DK}
\Br_A(A) \cong D_K(\Lambda)^{\iota,\mathrm{even}} / (1+\iota)D_K(\Lambda).
\end{equation}

\section{Complex tori with non-zero invariant Brauer group} \label{tori}

In this section we construct examples of simple complex tori of any dimension $g \geq 3$ with invariant Brauer group $\Z/2$.
We show that these complex tori are simple by calculating their endomorphism algebras.
We construct examples of two types:
\begin{enumerate}
\item[(i)] non-algebraisable complex tori with endomorphism ring $\Z$;
\item[(ii)] abelian varieties with endomorphism algebra an imaginary quadratic field.
\end{enumerate}
Note that it follows from Proposition~\ref{2.12} that an abelian variety with endomorphism ring $\Z$ always has trivial invariant Brauer group, since an abelian variety has $\rk_\Z(\NS\,(A)) \geq 1$.
Thus the examples (ii) have the smallest possible endomorphism algebras for abelian varieties with non-zero invariant Brauer group.

We use the following notation. The identity matrix of size $n$ is denoted by $I_n$.
We denote by $(1_{m},-1_n)$ the diagonal matrix with $m$ diagonal entries $1$
and $n$ diagonal entries $-1$, in this order. We write
$$\J=\left(\begin{array}{rr} 0& -1\\1&0\end{array}\right).$$
We use the direct sum symbol $\oplus$ for the direct sum of matrices. For example,
$\J^{\oplus n}$ is the square matrix of size $2n$ which is the direct sum of $n$ copies of $\J$.

\subsection{Complex tori with no non-trivial endomorphisms} \label{3.2}

The following theorem is the algebraic result which underlies our construction of non-algebraisable complex tori with endomorphism ring $\Z$ and invariant Brauer group~$\Z/2$.
The element $J$ is the complex structure giving rise to the desired torus, and the symmetric bilinear form $S$ is a generator for $\Bi_J(\Lambda)$ which we use in calculating the invariant Brauer group using Proposition~\ref{calc}.
The precursors of this theorem are \cite[Lemma 2]{Z91} and \cite[Lemma 2.5]{OZ95}.

\bthe \label{Y}
Let $g\geq 3$ and
let $\Lambda$ be a free abelian group of rank $2g$ with a non-degenerate symmetric bilinear form
$S\colon\Lambda\times\Lambda\to\Z$.
Let $\SO(\Lambda_\R,S)\subset\End(\Lambda_\R)$ be the special orthogonal group defined by $S$.
If the signature of $S$ is congruent to $2g\bmod 4$, 
then there is an element $J\in \SO(\Lambda_\R,S)$ such that $J^2=-1$
and the centraliser of $J$ in $\End(\Lambda_\Q)$ is $\Q\cdot \Id_\Lambda$.
\ethe
{\em Proof.} Let $\g=\mathfrak{o}(\Lambda_\R,S)$ be the Lie algebra of 
the real Lie group $G=\SO(\Lambda_\R,S)$.
Recall that
\[ \g = \{ A \in \End(\Lambda_\R) | S(Ax,y) + S(x,Ay) = 0 \text{ for all } x, y \in \Lambda_\R \}. \]
It is well known that $G$ can be given the structure of a {\em complete metric  space} 
such that
the group operations in $G$ are continuous.  For example, choosing a basis of the real vector space $\Lambda_{\R}$ we identify $G$ with a closed subset of the algebra of matrices 
$\mathrm{Mat}_{2g}(\R)\simeq\R^{(2g)^2}$ with its usual Euclidean metric. 

Let $C=\{J\in G|J^2=-1\}$. This is the set of all complex structures which preserve the symmetric bilinear form $S$. Our strategy is to show that, for very general $J \in C$, $\End_J(\Lambda_\Q)$ contains no elements other than $\Q \cdot \Id_\Lambda$.

The set $C$ is non-empty. Indeed, there is a basis of $\Lambda_\R$
with respect to which $S$ is given by the diagonal matrix 
$(1_m,-1_n)$. We have $m+n=2g$ and $m-n\equiv 2g\bmod 4$ (by the condition on signature in the statement of the theorem), hence $m$ and $n$ are both even.
Then the element of $\End(\Lambda_\R)$ whose matrix in this basis is $\J^{\oplus g}$
is in $C$, hence $C\not=\emptyset$.

We note that $C$ is invariant under
conjugations by the elements of $G$. We also note that $C\subset \g$, because for any $J\in C$ and
any $x,y\in \Lambda_\R$ we have
$$S(Jx,y)+S(x,Jy)=S(Jx,y)+S(J^{-1}x,y)=S(Jx,y)+S(-Jx,y)=0.$$
Since $C$ is a closed subset of $G$, it
inherits the structure of a complete metric space.
We shall show that, for each non-scalar endomorphism $u \in \End(\Lambda_\Q)$,
the closed subset $C_u=\{c\in C|cu=uc\}$ is nowhere dense in $C$.
This suffices to prove the theorem.
Indeed, by Baire's theorem a countable intersection of dense open
subsets of a complete metric space is dense.
Hence there exist elements $J \in C$ which are in the complement to all $C_u$ for non-scalar $u\in \End(\Lambda_\Q)$.
Any such $J$ commutes only with $\Q \cdot \Id_\Lambda$ in $\End(\Lambda_\Q)$, as required.

Thus let $u\in \End(\Lambda_\Q)$ be a non-scalar endomorphism. 
Let $\Ad\colon G\to \GL(\g)$ be the adjoint representation of the Lie group $G$ and let
$\ad\colon \g\to\gl(\g)$ be the adjoint representation of the Lie algebra $\g$.

Suppose for contradiction that $C_u$ has an interior point $c$, 
so there is an open neighbourhood $U$ of $c\in C$
such that $U\subset C_u$. Let $n$ be a positive integer. Recall that $C\subset \g$.
The function $G^n\to C$ sending $(g_1,\ldots,g_n)$ to 
$$\Ad(g_n)\ldots\Ad(g_1)c=g_n\ldots g_1cg_1^{-1}\ldots g_n^{-1}$$
is continuous, so there is an open neighbourhood $G_0$ of the unit element $e\in G$
such that $\Ad(g_n)\ldots\Ad(g_1)c\in U\subset C_u$ provided $g_i\in G_0$ for $i=1,\ldots,n$.
Hence $[\Ad(g_n)\ldots\Ad(g_1)c,u]=0$. This implies
$$[\ad(x_n)\ldots\ad(x_1)c,u]=[[x_n,[x_{n-1},\ldots[x_1,c]\ldots]],u]=0$$
for any $x_i\in\g$. This holds for any $n\geq 1$, hence the ideal of the Lie algebra $\g$ generated by $c$
commutes with $u$. 

Recall that $\g=\mathfrak{o}(\Lambda_\R,S)$ is an orthogonal Lie algebra of rank $2g$.
Since the rank is at least $6$, the Lie algebra $\mathfrak{o}(\Lambda_\R,S)$ is (absolutely) simple.
Therefore this ideal must be all of $\g$, thus $\g$ commutes with $u$.
However, the representation $\g\to\gl(\Lambda_\R)$ is absolutely irreducible, thus $u$ is scalar,
which is a contradiction. $\Box$

\smallskip

\bthe
For every integer $g \geq 3$, there exists a non-algebraisable complex torus $A$ of dimension $g$ with $\End(A) \cong \Z$ and $\Br_A(A) \cong \Z/2$.
\ethe

\noindent{\em Proof.}
Let $\Lambda$ be a free abelian group of rank $2g$.
Let $S$ be a non-degenerate, primitive, even symmetric bilinear form $S\colon\Lambda\times\Lambda\to\Z$
of signature congruent to $2g\bmod 4$, for example the form with matrix
$$ 2I_{2g-4} \oplus \left(\begin{array}{rr} 0& 1\\1&0\end{array}\right)^{\oplus 2}.$$
By Theorem \ref{Y}, if $g\geq 3$,
we can find a $J\in \SO(\Lambda_\R,S)$ such that $J^2=-1$ and the centraliser of $J$ in 
$\End(\Lambda_\Q)$ is $\Q\cdot{\rm Id}_\Lambda$.

Let $A$ be the complex torus $\Lambda_\R/\Lambda$, with complex structure $J$.
Thanks to (\ref{hom-tori}), $\End(A) = \End_J(\Lambda) \cong \Z$.

Since $J \in \SO(\Lambda_\R,S)$, we have $S \in \Bi_J(\Lambda)$.
By Lemma~\ref{bi-ends}, $\Bi_J(\Lambda) \otimes \Q = \Q S$.
Since $S$ is a primitive form, we conclude that $\Bi_J(\Lambda) = \Z S$.

Since $S$ is symmetric, $(1+\tau)\Bi_J(\Lambda) = 2\Z S$.
Since also $S$ is even, we conclude by Proposition~\ref{calc} that $\Br_A(A) = \Z/2$.

Note that $\Bi_J(\Lambda)$ contains no alternating forms, so $\NS(A) = \emptyset$ and hence $A$ is non-algebraisable.
$\Box$

\subsection{Simple abelian varieties with non-zero invariant Brauer group} \label{3.3}

The construction of simple abelian varieties with invariant Brauer group $\Z/2$ follows a more sophisticated version of the strategy from Section~\ref{3.2}.
We begin with an algebraic result.
In this theorem, $J$ is the complex structure associated with our desired abelian variety $A$, $S$ is an element of $\Bi_J(\Lambda)$ (indeed $S$ is ultimately a representative of the non-zero class in $\Br_A(A)$) and $J_0$ is a generator of $\End(A)$.

\bthe \label{H}
Let $g\geq 3$ be a positive integer.
Let $\Lambda_1$ and $\Lambda_2$ be free abelian groups of rank $2$ and $2g-2$,
respectively. Let $S_i\colon\Lambda_i\times\Lambda_i\to\Z$ be a positive-definite symmetric
bilinear form and let $J_i\in \SO(\Lambda_{i,\R},S_i)\cap\End(\Lambda_i)$
be such that $J_i^2=-\Id_{\Lambda_i}$, for $i=1,2$.
Let $S$ be the symmetric bilinear form on $\Lambda$ which is the orthogonal direct sum of $-S_1$ and $S_2$.
Let $J_0=J_1\oplus J_2\in \SO(\Lambda_\R,S)\cap\End(\Lambda)$.
Then there is an element $J\in \SO(\Lambda_\R,S)$ such that $J^2=-1$,
the centraliser of $J$ in $\End(\Lambda_\Q)$ is $\Q\,\Id_\Lambda+\Q J_0$,
and the symmetric bilinear form $S(Jx,J_0y)$ is positive-definite.
\ethe
{\em Proof.} Write $E(x,y)=S(x,J_0y)$. Then $E$ is an alternating form $\Lambda\times\Lambda\to\Z$ 
such that $E(x,y)=E(J_0x,J_0y)$.

Let $G$ be the centraliser of $J_0$ in $\SO(\Lambda_\R,S)$. We turn $\Lambda_\R$
into a complex vector space $\V$ via the complex structure $J_0$.
(This is not the complex structure we shall ultimately use to construct an abelian variety.)
Then $G\subset \GL(\V)$ is equal to the group of complex linear transformations 
preserving the Hermitian form $H(x,y)=S(x,y)+iE(x,y)$,
that is, the unitary group ${\rm U}(\V,H)$.
As before, the real Lie group $G$ can be given the structure of a complete metric space
such that the group operations in $G$ are continuous.

Let $\g={\mathfrak u}(\V,H)$ be the Lie algebra of the real Lie group $G$.
Recall that
$$\begin{array}{rcl} \g& =& \{ A \in \End_{\mathbb{C}}(\Lambda_\R) | H(Ax,y) + H(x,Ay) = 0 \text{ for all } x, y \in \Lambda_\R \} \\
&=&\{ A \in \End_{\mathbb{\R}}(\Lambda_\R) |  A J_0=J_0 A, \ S(Ax,y) + S(x,Ay) = 0 \text{ for all } x, y \in \Lambda_\R \}.
\end{array}$$
Define $C$ as the complement to $\{\pm J_0\}$ in
$$\{J\in \SO(\Lambda_\R,S)| J^2=-1, \ JJ_0 = J_0J\}=\{J\in {\rm U}(\V,H)| J^2=-1\}.$$
We exclude $\pm J_0$ here because they are isolated points in $\{J\in {\rm U}(\V,H)| J^2=-1\}$.
This implies that $C$ is a closed subset of $G$ and therefore is a complete metric space,
so we may use Baire's theorem as we did before.

For any $J\in C$ the bilinear form $S(Jx,J_0y)$ is symmetric. Thus the set
$$C^+=\{J\in C|S(Jx,J_0x)>0 \ \text{for all} \ x\not=0\}$$
is open in $C$. To show that $C^+\not=\emptyset$ we exhibit an element
of $C^+$. Namely, the linear transformation $J_0^\sharp=(-J_1) \oplus J_2$
preserves $S$, has determinant $1$ and commutes with $J_0=J_1\oplus J_2$. 
The symmetric form $S(J_0^\sharp x,J_0y)$
is the orthogonal direct sum of $S_1(x,y)$ and $S_2(x,y)$, hence is positive-definite.
Thus $J_0^\sharp\in C^+$, so $C^+$ is a non-empty open subset of $C$.

As before, $C$ is a subset of $\g$ invariant under conjugation by the elements of $G$.

For $u\in\End(\Lambda_\Q)$ let $C_u=\{c\in C| cu=uc\}$. We claim that if $C_u$ 
has an interior point $c$, then $u\in\Q\,\Id_\Lambda+\Q J_0$. Indeed, 
arguing as in the proof of Theorem \ref{Y} we obtain that
the ideal of $\g$ generated by $c$ commutes with $u$.
The Lie algebra $\g$ is the direct sum of its centre $Z(\g)$ and the simple Lie algebra
$\mathfrak{su}(\V,S)$. We have $Z(\g)=i\R\Id_\V=\R J_0$.
The only $J \in \R J_0$ which satisfy $J^2 = -1$ are $\pm J_0$, which are not in $C$, so $c\notin Z(\g)$.
Thus the ideal of $\g$ generated by $c$ contains $\mathfrak{su}(\V,S)$.
A direct calculation shows that the centraliser of $\mathfrak{su}(\V,S)$ in $\End(\Lambda_\R)=
\gl(\V)\oplus\gl(\V)\tau$, where $\tau$ is the complex conjugation,
is $\C\,\Id_\V\subset \gl(\V)$. (For this calculation it is crucial that $g\geq 3$.) 
Thus $u\in\Q\,\Id_\Lambda+\Q J_0$. 

Let $C_0$ be the union of $C_u$ for all $u\in\End(\Lambda_\Q)$ such that 
$u\notin\Q\,\Id_\Lambda+\Q J_0$.
By Baire's theorem, the complement to $C_0$ is dense in $C$. Hence we can find an element
$J\in C^+$ outside of $C_0$. It has all the required properties. $\Box$

\smallskip

\bthe
For every integer $g \geq 3$, there exists a complex abelian variety $A$ of dimension $g$ with $\End(A) \otimes \Q \cong \Q(i)$ and $\Br_A(A) \cong \Z/2$.
\ethe

\noindent {\em Proof.}
Let $M$ be a symmetric, positive definite $2\times 2$-matrix with entries in $\Z$
such that the diagonal entries are even. Assume that $M\not=nM'$, where $n\geq 2$ is an integer
and $M'$ has entries in $\Z$. For example, we can take
$$M=\left(\begin{array}{rr} 2& 1\\1&2\end{array}\right).$$ 

Let $\Lambda$ be a free abelian group 
which is a direct sum $\Lambda_1\oplus\Lambda_2$, where $\rk(\Lambda_1)=2$ and 
$\rk(\Lambda_2)=2g-2$.
Choose a $\Z$-basis in each $\Lambda_i$. Let $S_1$ be the symmetric bilinear
form on $\Lambda_1$ with matrix $2I_2$. Let $S_2$ be the symmetric bilinear form on $\Lambda_2$
with matrix $M^{\oplus 2}\oplus 2I_{2g-6}$. It is clear that the orthogonal direct sum 
$S=(-S_1)\oplus S_2$ is non-degenerate, even and primitive. Let
$$J_0=\J\oplus\left(\begin{array}{cc} 0& -I_2\\I_2&0\end{array}\right)\oplus \J^{g-3}.$$
This is chosen to have the following properties: $J_0$ preserves the decomposition $\Lambda=\Lambda_1\oplus\Lambda_2$
and satisfies $J_0^2=-1$ and $S(J_0x,J_0y)=S(x,y)$ for all $x,y\in \Lambda$.
Hence by Theorem \ref{H} there is an element
$J\in \SO(\Lambda_\R,S)$ such that $J^2=-1$,
the centraliser of $J$ in $\End(\Lambda_\Q)$ is $\Q\,\Id_\Lambda+\Q J_0$,
and the symmetric form $S(Jx,J_0y)$ is positive-definite.

Let $A$ be the complex torus $\Lambda_\R/\Lambda$ associated with the complex structure $J$ on~$\Lambda_\R$.
By (\ref{hom-tori}), $\End(A) \otimes \Q = \End_J(\Lambda_\Q) = \Q\,\Id_\Lambda \oplus \Q J_0 \cong \Q(i)$.

Write $E(x,y)=S(x,J_0y)$. The matrix of $E$ is 
$$(-2\J)\oplus\left(\begin{array}{cc} 0& -M\\M&0\end{array}\right)\oplus (2\J)^{g-3}.$$

We have $S, E \in \Bi_J(\Lambda)$.
Thanks to Lemma~\ref{bi-ends}, $\Bi_J(\Lambda) \otimes \Q = \Q S + \Q E$.
Looking at the matrices of $S$ and $E$
we observe that the set of matrices in $\Q S+\Q E$ that have all entries in $\Z$ is precisely
$\Z S+\Z E$. Hence $\Bi_J(\Lambda)=\Z S+\Z E$.

The form $S$ is symmetric while $E$ is anti-symmetric.
Hence $(1+\tau)\Bi_J(\Lambda) = 2\Z S$.
Since $S$ is even, we conclude from Proposition~\ref{calc} that $\Br_A(A) \cong \Z_2 S/2\Z_2 S \cong \Z/2$.

Finally $E(Jx,y)=S(Jx,J_0y)$ is positive definite.
Hence by (\ref{polarisation}), $E \in \Bi_J(\Lambda)$ corresponds to a polarisation of $A$ and so $A$ is an abelian variety.
$\Box$

\section{Complex multiplication} 

In this section we study the invariant Brauer group of complex abelian varieties of CM type.
We construct two different examples of abelian surfaces of CM type with invariant Brauer group $\Z/2$: one is isogenous to a product of elliptic curves, the other is simple.

\subsection{Non-simple abelian surfaces of CM type} \label{surfaces}

Let us first give a construction of abelian surfaces $A$ isogenous to the product of two elliptic curves with CM
and $\Br_A(A) \cong \Z/2$. The subtlety here is that, on the one hand, 
the two imaginary quadratic fields
cannot be the same (Corollary \ref{c2}) and, on the other hand, $A$ cannot be
a product of two elliptic curves (Corollary~\ref{c1}).

Let $d_1$ and $d_2$ be square-free negative integers such that 
$d_1\equiv d_2\equiv -1\bmod 4$. Let 
$$\Lambda_i=\Z[\sqrt{d_i}]=\{x+y\sqrt{d_i}|x,y\in\Z\},
\quad\text{for} \quad i=1,2.$$
Equip $\Lambda_{i,\R}$ with the obvious complex structure.
Then $E_i=\Lambda_{i,\R}/\Lambda_i$ is an elliptic curve over $\C$
such that $\End(E_i)=\Lambda_i$.

Since $2$ is ramified in $\Q(\sqrt{d_i})$, we have
$$\Lambda_i/2=\F_2[\sqrt{d_i}]\cong\F_2[x]/(x^2)$$ for $i=1,2$, 
where the second isomorphism identifies $x$ with $1+\sqrt{d_i}$.
Let $p_1, p_2$ denote the composed homomorphisms
$\Lambda_i \to \Lambda_i/2 \to \F_2[x]/(x^2)$ and let
$$ \Lambda = \{ (x_1, x_2) \in \Lambda_1\oplus\Lambda_2 | p_1(x_1) = p_2(x_2) \}. $$
Then $\Lambda$ is a subring of $\Lambda_1 \oplus \Lambda_2$ and, as a lattice, $\Lambda$ is generated by $(2,0)$, $(1,1)$, $(2\sqrt{d_1}, 0)$ and $(\sqrt{d_1}, \sqrt{d_2})$.

Let $A=\Lambda_\R/\Lambda$, where $\Lambda_\R$ is equipped with the complex structure $J$ which is the direct sum of the complex structures on $\Lambda_{1,\R}$ and $\Lambda_{2,\R}$. Thus
$A=(E_1\times E_2)/G$, where $G$ is the graph of 
the isomorphism of abelian groups $E_1[2]\tilde\lra E_2[2]$
which is the same as $\Lambda_1/2 \tilde\lra \F_2[x]/(x^2) \tilde\lra \Lambda_2/2$.

Let $K=\Q(\sqrt{d_1})\oplus \Q(\sqrt{d_2})$.
If $d_1 \neq d_2$, then the elliptic curves $E_1$ and $E_2$ are not isogenous so $\End(A) \otimes \Q \cong K$. 
Since $\Lambda$ is a conjugation-invariant subring of $K$, 
we can use \eqref{DK-subring} to calculate
\[ D_K(\Lambda) = (\tfrac{1}{2}, 0)\Z + (\tfrac{1}{4}, -\tfrac{1}{4})\Z + (\tfrac{1}{2\sqrt{d_1}}, 0)\Z + (\tfrac{1}{4\sqrt{d_1}}, \tfrac{-1}{4\sqrt{d_2}})\Z. \]

Now $D_K(\Lambda)^\iota$ has a basis $e_1 = (\tfrac{1}{2}, 0)$ and $e_2 = (\frac{1}{4}, -\tfrac{1}{4})$.
For $a, b \in \Q$, we can calculate
\begin{align*}
B_{ae_1 + be_2}((2,0), (2,0)) &= 4a+2b,
\\ B_{ae_1 + be_2}((1,1), (1,1)) &= a,
\\ B_{ae_1 + be_2}((2\sqrt{d_1}, 0), (2\sqrt{d_1}, 0)) &= -4d_1a-2d_1b,
\\ B_{ae_1 + be_2}((\sqrt{d_1}, \sqrt{d_2}), (\sqrt{d_1}, \sqrt{d_2})) &= -d_1a + \tfrac{1}{2}(d_2 - d_1)b.
\end{align*}
Recalling that $d_1 \equiv d_2 \bmod 4$, we conclude that $B_{e_1}$ is odd while $B_{e_2}$ is even.
Hence $D_K(\Lambda)^{\iota,\mathrm{even}} = 2\Z e_1 + \Z e_2$.
Meanwhile $(1+\iota)D_K(\Lambda) = 2\Z e_1 + 2\Z e_2$.
Hence by \eqref{DK},
$\Br_A(A) = \Z/2$, generated by $e_2$.

\brem \label{ppas} {\rm
The bilinear form $B_{(1/4\sqrt{d_1}, 1/4\sqrt{d_2})} \in \Bi_J(\Lambda)$ is alternating and unimodular and the associated Hermitian form is positive definite.
Hence $A$ has a principal polarisation.
This shows that there are {\em principally polarised} abelian surfaces with non-trivial invariant Brauer group.}
\erem

\subsection{Simple abelian surfaces of CM type} \label{4.2}

We now give examples of simple abelian surfaces $A$ of CM type with $\Br_A(A) \cong \Z/2$.

Let $K$ be a quartic CM field which is not biquadratic over $\Q$, but such that $K_2 = K \otimes \Q_2$ is a biquadratic extension of $\Q_2$.
(See below for an example of such a field.)
Since $K$ is a quartic CM field, it contains a real quadratic field $K_+ = \Q(\sqrt{d})$ where $d \in \Q$.
Since $K_2/\Q_2$ is biquadratic, we can write $K_2 = \Q_2(\sqrt{d}, \sqrt{e})$ for some $e \in \Q_2$ (where $e \not\in \Q^\times K_+^{\times 2}$ since $K/\Q$ is not biquadratic, and $e, e/d \not\in \Q_2^{\times 2}$ since $\sqrt{e} \not\in \Q_2(\sqrt{d})$).

Let $\Lambda$ be a $\Z$-lattice in $K$ such that
\[ \Lambda_2 = \Lambda \otimes_\Z \Z_2 = \Z_2 + \Z_2\sqrt{d} + \Z_2\sqrt{e} + \Z_2\sqrt{de}. \]
Let $A = K_\R / \Lambda$.
Choosing $J \in K_\R$ such that $J^2 = -1 $ gives $A$ the structure of a complex torus.
The argument of \cite[pp.~212--213]{M} shows that $A$ is polarisable and hence an abelian variety.
(Unlike in \cite[pp.~212--213]{M}, we do not have $\Lambda = \O_K$.
However $n\Lambda \subset \O_K$ for some positive integer~$n$ so we can replace the value $\alpha$ from \cite[p.~212]{M} by $n^2 \alpha$ to ensure that $\tr_{K/\Q}(\alpha x \bar y) \in \Z$ for all $x,y \in \Lambda$.)

Since $K/\Q$ is not biquadratic, $K_+$ is the only quadratic subfield of $K$.
Hence $K$ does not contain any imaginary quadratic fields, so every CM type for $K$ is primitive.
Hence $A$ is a simple abelian variety with $\End(A) \otimes \Q = K$,
by Shimura and Taniyama, see \cite[Thm.~3.5]{Lang}.

Write $D_K(\Lambda_2) = D_K(\Lambda) \otimes \Z_2$.
Since $\Lambda_2$ is a subring of $K_2$, $\alpha \in D_K(\Lambda_2)$ if and only if $\tr_{K_2/\Q_2}(\alpha x) \in \Z_2$ for all $x \in \Lambda_2$, from which we obtain
\[ D_K(\Lambda_2) = \frac{1}{4}\Z_2 + \frac{1}{4\sqrt{d}}\Z_2 + \frac{1}{4\sqrt{e}}\Z_2 + \frac{1}{4\sqrt{de}}\Z_2. \]

If $\alpha = a_1 + a_2\sqrt{d} + a_3\sqrt{e} + a_4\sqrt{de}$ ($a_1, a_2, a_3, a_4 \in \Q_2$), then
\begin{gather*}
B_\alpha(1,1) = 4a_1, \quad B_\alpha(\sqrt{d}, \sqrt{d}) = 4da_1,
\\ B_\alpha(\sqrt{e}, \sqrt{e}) = -4ea_1, \quad B_\alpha(\sqrt{de}, \sqrt{de}) = -4dea_1.
\end{gather*}
Hence $B_\alpha$ is even if and only if $4a_1 \in 2\Z_2$ or in other words, $a_1 \in \frac{1}{2}\Z_2$.
Consequently
\[ D_K(\Lambda_2)^{\iota,\mathrm{even}} = \frac{1}{2}\Z_2 + \frac{1}{4\sqrt{d}}\Z_2 \]
and
\[ (1+\iota) D_K(\Lambda_2) = \frac{1}{2}\Z_2 + \frac{1}{2\sqrt{d}}\Z_2. \]
Hence by \eqref{DK}, $\Br_A(A) = \Z/2$, generated by $1/4\sqrt{d}$.

\ble
Let $K_+ = \Q(\sqrt{5})$ and $K = K_+(\sqrt{\delta})$ where $\delta = -30 + 8\sqrt{5}$.
Then $K$ is a quartic CM field, $K/\Q$ is not biquadratic, but $K \otimes \Q_2/\Q_2$ is biquadratic.
\ele

\textit{Proof.}
We have $-30 \pm 8\sqrt{5} < 0$ so $\delta$ is a totally negative element of $K_+$.
Hence $K$ is a CM field.

If $K/\Q$ were biquadratic, then $\delta \in \Q^\times K_+^{\times 2}$ so $\mathrm{Nm}_{K_+/\Q}(\delta) \in \Q^{\times 2}$.
But $\mathrm{Nm}_{K_+/\Q}(\delta) = 30^2 - 8^2 \times 5 = 580$ which is not a square.

Now we examine $K \otimes \Q_2$.
First note that $K_+ \otimes \Q_2 = \Q_2(\sqrt{5})$ is the unramified quadratic extension of $\Q_2$.
We have $\delta = 2(-15 + 4\sqrt{5})$ so the $2$-adic valuation of $\delta$ is odd.
Hence $\delta$ does not have a square root in $K_+ \otimes \Q_2$, so $K \otimes \Q_2$ is a field.

Finally $(2+\sqrt{5})^2 = 9+4\sqrt{5}$.
We have $-15 + 4\sqrt{5} \equiv 9 + 4\sqrt{5} \bmod 8$ so by Hensel's lemma, $-15 + 4\sqrt{5} = \zeta^2$ for some $\zeta \in \Q_2(\sqrt{5})$.
Therefore
\[ K \otimes \Q_2 = \Q_2(\sqrt{5})(\sqrt{2\zeta^2}) = \Q_2(\sqrt{5})(\sqrt{2}) \]
is a biquadratic extension of $\Q_2$.
$\Box$

\subsection{Simple abelian varieties of CM type} \label{simple-cm}

In this section we show that any abelian variety of CM type whose endomorphism ring is the ring of integers of a cyclotomic field has $\Br_A(A) = 0$.

Let $K$ be a CM field. Write $[K:\Q]=2g$ so that there is an isomorphism of $\R$-algebras
$K_\R=K\otimes_\Q\R\cong\C^g$. Let $\O_K$ be the ring of integers of $K$.
Then $A=K_\R/\O_K$ is a real torus.

Suppose that $J\in K_\R$ is such that $J^2=-1$ and $J$ is not contained in a proper CM subfield of $K$
(that is, the associated CM type is {\em primitive}). The action of $J$ makes $K_\R$ a complex
vector space so $A$ is a complex torus such that $\End(A)=\O_K$. It is well known that $A$ being a simple
abelian variety is equivalent to the CM type being primitive \cite[Ch.~22]{M}.

Let $\D_K\subset\O_K$ be the {\em different}.
Recall that the fractional ideal $\D^{-1}_K\subset K$ is the dual of $\O_K$ with respect
to the bilinear form $\tr_{K/\Q}(xy)$.
Hence by \eqref{DK-subring},
\[ D_K(\O_K) = \D_K^{-1}. \]
Write \( (\D_K^{-1})_+ \) for the subgroup of all conjugation-invariant elements of \( \D_K^{-1} \).

Let us consider the particular case when $K=\Q(\zeta_n)$, where $\zeta_n$ is a
primitive root of unity of degree $n\geq 3$. Since $\O_k=\Z[\zeta_n]$,
we see that $\{\zeta_n^i|i=0,\ldots,\varphi(n)-1\}$ is a $\Z$-basis of $\O_K$.
For any $x=\zeta_n^i$ we have $B_\alpha(x,x)=\tr_{K/\Q}(\alpha x\bar x)=\tr_{K/\Q}(\alpha)$.
Thus $B_\alpha(x,y)$ is an even form if and only if $\tr_{K/\Q}(\alpha)$ is even. So
in this case, from \eqref{DK}, we have
$$\Br_A(A)=\{\alpha\in (\D_K^{-1})_+|\ \tr_{K/\Q}(\alpha)\equiv 0\bmod 2\}/(1+\iota)\D_K^{-1}.$$

\ble \label{odin} Let $K=\Q(\zeta_n)$, where $n\geq 3$.
If $n$ is odd, then the different ideal of $K$ is 
the principal ideal of $\O_K$ generated by 
$$n\prod_{\text{\rm primes}\, p|n}\frac{1}{\zeta_p-\zeta_p^{-1}}.$$
If $n$ is even, then
the different ideal of $K$ is generated by 
$$\frac{n}{2}\prod_{\text{\rm odd primes} \, p|n}\frac{1}{\zeta_p-\zeta_p^{-1}}.$$
\ele
{\em Proof.} If $p$ is a prime and $r$ is a positive integer, then the different of 
$\Q(\zeta_{p^r})$ is the ideal of $\O_{\Q(\zeta_{p^r})}$ generated by 
$(1-\zeta_{p^r})^{p^{r-1}((p-1)r-1)}$, see \cite[Ch. 3, Thm.~27]{SD}. Using that 
$p\O_{\Q(\zeta_{p^r})}=(1-\zeta_{p^r})^{p^{r-1}(p-1)}\O_{\Q(\zeta_{p^r})}$
and $(1-\zeta_p)\O_{\Q(\zeta_{p^r})}=(1-\zeta_{p^r})^{p^{r-1}}\O_{\Q(\zeta_{p^r})}$
we obtain that the different is generated by $p^r(1-\zeta_p)^{-1}$.

If $p$ is odd, then $1+\zeta_p$ is a unit in $\O_{\Q(\zeta_{p^r})}$ so the different is generated by
$$ p^r(1-\zeta_p)^{-1}(1+\zeta_p)^{-1} = \frac{p^r}{\zeta_p-\zeta_p^{-1}}. $$
If $p=2$, then $p^r(1-\zeta_p)^{-1} = p^r/2$.
Thus the lemma holds whenever $n$ is equal to a prime power $p^r$.

If $m$ and $m'$ are coprime positive integers, then the fields $\Q(\zeta_m)$ and $\Q(\zeta_{m'})$
are linearly disjoint with coprime discriminants, hence by \cite[III.2.13, VI.1.14]{FT} we have
$\O_{\Q(\zeta_{mm'})}=\O_{\Q(\zeta_m)}\otimes \O_{\Q(\zeta_{m'})}$. 
This implies that the different of $K$ is the product of the differents of the subfields 
$\Q(\zeta_{p^r})\subset K$,
where $r$ is the highest power of $p$ dividing $n$. $\Box$

\medskip

The following lemma is a version of \cite[Sec.~3.3, example~5]{Ser77} over $\Z$.

\ble \label{ind-mods}
Let $G$ be a group, let $M$ be a $\Z[G]$-module which is free as a $\Z$-module, and let $N$ be a free $\Z[G]$-module.  Then $M \otimes_\Z N$ is a free $\Z[G]$-module.
\ele
{\em Proof.}
Let $\{m_i | i\in I\}$ be a $\Z$-basis for $M$, and let $\{n_j | j \in J\}$ be a $\Z[G]$-basis for~$N$.
For each $g \in G$, $\{gm_i | i \in I\}$ is also a $\Z$-basis for $M$.  Consequently for each $g \in G$ and $j \in J$, $\{ gm_i \otimes gn_j | i \in I \}$ is a $\Z$-basis for $M \otimes_\Z gn_j\Z$.
Since $\{gn_j | g \in G, j \in J\}$ is a $\Z$-basis for $N$, we deduce that $$\{ gm_i \otimes gn_j | g \in G, i \in I, j \in J \}$$ is a $\Z$-basis for $M \otimes_\Z N$.
This basis is $G$-stable so $M \otimes_\Z N$ is free as a $\Z[G]$-module.
$\Box$

\ble \label{dva}
Let $n\geq 3$ be an integer which is not a power of $2$. Then $\O_{\Q(\zeta_{n})}$ is a 
free $\Z[\Z/2]$-module, where $\Z/2$ acts by complex conjugation.
\ele
{\em Proof.}
Let $p$ be an odd prime and let $r$ be a positive integer.
Let 
\begin{align*}
   S     &=\{\zeta_{p^r}^j|1\leq j\leq p^{r-1}(p-1)/2\},
\\ \ov S &=\{\zeta_{p^r}^j|p^{r-1}(p+1)/2 \leq j \leq p^r-1\}.
\end{align*}
We claim that $S\cup\ov S$ is a basis of $\O_{\Q(\zeta_{p^r})}=\Z[\zeta_{p^r}]$.
Indeed, $|S\cup\ov S|=\varphi(p^r)$, so it is enough to show that every power of $\zeta_{p^r}$
is in the subgroup $N$ of $\Z[\zeta_{p^r}]$ generated by
$S\cup\ov S$. The set $S\cup\ov S$ contains $\zeta_{p^r}^{p^{r-1}a}$, where $a=1,\ldots, p-1$, 
hence $1\in N$. For $j=1,\ldots, p^r-1$ the relation $\sum_{a=0}^{p-1}\zeta_{p^r}^{j+p^{r-1}a}=0$
shows that each power of $\zeta_{p^r}$ between $p^{r-1}(p-1)/2$ and $p^{r-1}(p+1)/2$
is also in $N$, so we are done.
This proves the statement for $n=p^r$.

Now let $n\geq 3$ be an integer which is not a power of 2.
Then we can write $n = p^r m$ where $p$ is an odd prime, $r$ is a positive integer, and $p \nmid m$.
The fields $\Q(\zeta_{p^r})$ and $\Q(\zeta_m)$
are linearly disjoint with coprime discriminants, hence 
$\O_{\Q(\zeta_n)}=\O_{\Q(\zeta_{p^r})}\otimes \O_{\Q(\zeta_m)}$. 
We have proved that $\O_{\Q(\zeta_{p^r})}$ is a free $\Z[\Z/2]$-module, so we can conclude by applying Lemma~\ref{ind-mods}. $\Box$

\bpr \label{cyclo}
Let $K=\Q(\zeta_n)$, where $n\geq 3$.
Let $A$ be an abelian variety attached to $K$ as at the beginning of Section~\ref{simple-cm}, that is, $A = K_\R/\O_K$ with a primitive CM type.
Then $\Br_A(A)=0$.
\epr
{\em Proof.} Let us assume that $n$ is not a power of $2$. By Lemma \ref{odin} we have
 $\D_K^{-1}=\eta\O_K$ with $\eta$ a rational multiple of $\prod_p(\zeta_p-\zeta_p^{-1})$,
where $p$ ranges over the odd prime factors of $n$. We note that $\bar\eta=\pm\eta$.
If $\bar\eta=\eta$, then $\Br_A(A)$ is
$$(\D_K^{-1})^{\rm even}_+/(1+\iota)\D_K^{-1}\subset
(\D_K^{-1})_+/(1+\iota)\D_K^{-1}\cong(\O_K)_+/(1+\iota)\O_K=\widehat\H^0(\Z/2,\O_K).$$
By Lemma \ref{dva} the $\Z[\Z/2]$-module $\O_K$ is free, hence $\widehat\H^0(\Z/2,\O_K)=0$.

If $\bar\eta=-\eta$, then $\Br_A(A)$ is
$$(\D_K^{-1})^{\rm even}_+/(1+\iota)\D_K^{-1}\subset
(\D_K^{-1})_+/(1+\iota)\D_K^{-1}\cong(\O_K)_-/(1-\iota)\O_K=\H^1(\Z/2,\O_K).$$
This is also zero by Lemma \ref{dva}.

Now let $n=2^m$, where $m\geq 2$. In this case $\D_K^{-1}$ is generated by $2^{1-m}$,
hence
$$(\D_K^{-1})_+/(1+\iota)\D_K^{-1}\cong(\O_K)_+/(1+\iota)\O_K=\widehat\H^0(\Z/2,\O_K).$$
Then
$\O_K=\Z[\zeta_{2^m}]$ has a $\Z$-basis 
$\{1\}\cup\{i\}\cup R\cup\ov R$, where $R=\{\zeta_{2^m}^j| 1\leq j\leq 2^{m-2}-1\}$.
The submodule spanned by $R\cup\ov R$ is a free $\Z[\Z/2]$-module so 
$\widehat\H^0(\Z/2,\Z R \oplus \Z \ov R) = 0$.
Furthermore, we have $\widehat\H^0(\Z/2,\Z i) = 0$.
It follows that $\widehat\H^0(\Z/2,\O_K)\cong\Z/2$
generated by the class of $1$. Thus $(\D_K^{-1})_+/(1+\iota)\D_K^{-1}$ is generated by the 
class of $2^{1-m}$. We note that $\tr_{K/\Q}(2^{1-m})=1$, so $2^{1-m}$ is not contained in
$(\D_K^{-1})^{\rm even}$. Therefore,
$\Br_A(A)=(\D_K^{-1})^{\rm even}_+/(1+\iota)\D_K^{-1}=0$. $\Box$

\bthe
Let $n \ge 3$ be an integer, let $K=\Q(\zeta_n)$ be the $n$-th cyclotomic field and let $\O_K$ be the ring of integers of $K$. Let $\varphi(n):=[K:\Q]$
and $g=\varphi(n)/2$.
If $A$ is a simple $g$-dimensional abelian variety of CM type over $\C$ with multiplication by $\O_K$, then
$\Br_A(A)=0$.
\ethe
{\em Proof.} We may assume that $A(\C)=K_{\R}/I$ where $I=\mathrm{H}_1(A(\C),\Z)$ is a non-zero ideal in $\O_K$.
We claim that there is a non-zero {\em principal} ideal $J \subset I$ such that the (finite) quotient $I/J$ has {\em odd} order. Then there is a simple complex abelian variety 
$B$ of CM type with multiplication by $\O_K$ such that 
$J=\mathrm{H}_1(B(\C),\Z)$, $B(\C)=K_{\R}/J$ and there is an isogeny $B \to A$ with kernel $I/J$. In particular, the degree
of $B \to A$ is odd.
Since $A$ is simple, so is $B$, and thus $B$ falls within the scope of Proposition~\ref{cyclo}.
Hence $\Br_B(B) = 0$, and it follows from Proposition \ref{isogOdd} that $\Br_A(A)=0$.

Let us prove the existence of such an ideal $J$. According to \cite[Thm.~18.20]{CR62}, there exists a non-zero ideal $M\subset \O_K$ such that
$M$ is coprime to the ideal $2\O_K$ and $J:=IM$ is a principal ideal in $\O_K$. Clearly, $J\subset I$.  On the other hand, according to \cite[Cor.~18.24]{CR62},
the additive groups $\O_K/M$ and $I/IM=I/J$ are isomorphic. In particular, the order of $O_K/M\cong I/J$  is odd. $\Box$

\section{Abelian varieties over non-closed fields} \label{non-closed}

Let $A$ be an abelian variety over a field $k$. The Picard scheme 
${\bf Pic}_{A/k}$ is a group $k$-scheme which is smooth at $e$
by \cite[Ch.~VI, Thm.~6.18]{vdGM} and hence is smooth over $k$. Since $A$ has a $k$-point,
${\bf Pic}_{A/k}$ represents the relative Picard functor $\Pic_{A/k}$ that sends
a $k$-scheme $T$ to $\Pic(A\times_kT)/\Pic(T)$, see
\cite[Thm.~2.5]{kleiman}.
Thus for any field extension $k\subset L$ we have 
${\bf Pic}_{A/k}(L)=\Pic(A_L)$, where $A_L=A\times_kL$.

Let ${\bf Pic}^0_{A/k}$ be the connected component of $e$ in ${\bf Pic}_{A/k}$.
This is the dual abelian variety $A^\vee$, see \cite[Ch.~VI, (6.19)]{vdGM}.
Write $\Pic^0(A_L)={\bf Pic}^0_{A/k}(L)$.
The N\'eron--Severi group of $A_L$ is defined as $\NS(A_L)=\Pic(A_L)/\Pic^0(A_L)$.
In particular, $\NS(\ov A)$ is the group of connected components of ${\bf Pic}_{A/k}\times_k\bar k$.

The group $k$-scheme of connected components of ${\bf Pic}_{A/k}$ is \'etale
\cite[VI$_{A}$, 5.5]{SGA3}. Thus the connected components of 
${\bf Pic}_{A/k}\times_k\bar k$ are obtained from the connected components of
${\bf Pic}_{A/k}\times_k k_\s$ by base change from $k_\s$ to $\bar k$.
The smoothness of the $k_\s$-scheme ${\bf Pic}_{A/k}\times_kk_\s$ 
implies that each connected component of 
${\bf Pic}_{A/k}\times_kk_\s$ is smooth over $k_\s$, hence contains a $k_\s$-point.
It follows that the natural map $\NS(A_\s)\to \NS(\ov A)$ is an isomorphism.
See also \cite[Ch.~III, (3.27), (3.29)]{vdGM}.

By a theorem of Chow, for abelian varieties $A$ and $B$ over $k$
the natural map $\Hom(A_\s,B_\s)\to \Hom(\ov A,\ov B)$ is an isomorphism
(see  \cite[Thm.~3.19]{Con06} for a modern proof).

Recall that a line bundle $L$ on $\ov A$ defines a canonical map
$\varphi_L\colon\ov A\to \ov A^\vee$, see \cite[Ch.~6, p.~60]{M}.
Let $\mathcal P$ be the Poincar\'e line bundle on $A\times A^\vee$.
We write $\pi_j\colon A\times A\to A$ for the projection to the $j$-th factor, $j=1,2$.
Let $i_1\colon A\to A\times A$ be the map sending $x$ to $(x,e)$, and similarly $i_2(x)=(e,x)$.

The following lemma is essentially well known.

\ble \label{known}
{\rm (i)} There is a commutative diagram of abelian groups
$$\begin{array}{ccccccc}
\Pic(\ov A)&\oplus &\Pic(\ov A)&\oplus &\Hom(\ov A,\ov A^\vee)&\stackrel{\sim}\lra&
\Pic(\ov A\times\ov A)\\
\uparrow&&\uparrow&&||&&\uparrow\\
\Pic(A_\s)&\oplus& \Pic(A_\s)&\oplus& \Hom(A_\s,A^\vee_\s)&\stackrel{\sim}\lra&\Pic(A_\s\times A_\s)
\end{array}$$
where in each row the first two summands are embedded via $\pi_1^*$ and $\pi_2^*$,
the third summand is embedded via the map sending $\phi\in \Hom(\ov A,\ov A^\vee)$
to $({\rm id},\phi^*)\mathcal P$, and the vertical arrows are injective.

{\rm (ii)} There is a commutative diagram of abelian groups
$$\begin{array}{ccccccccc}
\Pic(\ov A)&\twoheadrightarrow&\NS(\ov A)&\stackrel{\sim}\lra& 
\Hom(\ov A,\ov A^\vee)^{\rm sym}&\hookrightarrow&
\Hom(\ov A,\ov A^\vee)&\hookrightarrow&\Pic(\ov A\times\ov A)\\
\uparrow&&||&&||&&||&&\uparrow\\
\Pic(A_\s)&\twoheadrightarrow&\NS(A_\s)&\stackrel{\sim}\lra& 
\Hom(A_\s,A_\s^\vee)^{\rm sym}&\hookrightarrow&
\Hom(A_\s,A_\s^\vee)&\hookrightarrow&\Pic(A_\s\times A_\s)
\end{array}
$$
where in each row the composed map is $\delta=m^*-\pi_1^*-\pi_2^*$ and the second arrow 
is $L\mapsto\varphi_L$.
\ele
{\em Proof.} (i) Injectivity of the vertical arrows follows from the fact that ${\bf Pic}_{A/k}$
represents the relative Picard functor. 
The fact that the horizontal maps are isomorphisms is a consequence of the obvious
properties  of the maps $\pi_1$, $\pi_2$, $i_1$, $i_2$, and the
universal property of the Poincar\'e line bundle, cf. \cite[Prop.~1.7]{SZ14}.

(ii) Let $L\in\Pic(\ov A)$. If $L\in\Pic^0(\ov A)$, then
$\delta(L)=m^*L\otimes p_1^*L^{-1}\otimes p_2^*L^{-1}=0$ by the basic theory of
abelian varieties, so $\delta$ factors through $\NS(\ov A)$.
Since $i_1^*\delta=i_2^*\delta=0$, we see that
$\delta$ maps $\Pic(\ov A)$ to $\Ker(i_1^*)\cap\Ker(i_2^*)$, which by part (i) is identified with
$\Hom(\ov A,\ov A^\vee)$ by the map sending $\phi\in \Hom(\ov A,\ov A^\vee)$
to $({\rm id},\phi^*)\mathcal P$. By \cite[Ch.~8]{M} there is
a canonical isomorphism
$$({\rm id},\varphi_L^*)\mathcal P=m^*L\otimes p_1^*L^{-1}\otimes p_2^*L^{-1}.$$ 
The same arguments apply to the bottom row. $\Box$

\medskip

Recall that $e\colon\Spec(k)\to A$ 
denotes the unit element of $A$. In Section \ref{defi} we defined
$\Br_e(A)=\Ker[e^*\colon\Br(A)\to\Br(k)]$. Let 
$$\Br_{a}(A):=\Br_e(A)\cap\Br_{1}(A), \quad\quad \Br_{a,A}(A):=\Br_e(A)\cap\Br_{1}(A)\cap\Br_A(A).$$
Let $\bar k$ be an algebraic  closure of $k$, and let $k_\s$ be the separable closure of $k$ in $\bar k$,
with Galois group $\Ga=\Gal(k_\s/k)$. Write
$\ov A=A\times_k\bar k$ and $A_\s=A\times_k k_\s$. The Leray spectral sequence
\cite[Thm.~III.1.18 (a)]{EC}
\begin{equation}
\H^i(k,\H_\et^j(A_\s,\G_m))\Rightarrow\H_\et^{i+j}(A,\G_m) \label{ss}
\end{equation}
gives rise to an exact sequence
$$\Br(k)\lra\Br_1(A)\lra \H^1(k,\Pic(A_\s))\lra \H^3(k,\G_m)\lra \H_\et^3(A,\G_m).$$
The first arrow composed with $e^*$ is the identity map on $\Br(k)$.
Similarly, the composition of the last arrow with the map induced by $e$ is the identity map
on $\H^3(k,\G_m)$. Thus we obtain canonical isomorphisms
$$\Br_{a}(A)\cong\H^1(k,\Pic(A_\s))\cong\Br_1(A)/\Br_0(A).$$

\bpr \label{rain}
Let $A$ be an abelian variety over a field $k$. Then $\Br_{a,A}(A)$ is canonically isomorphic
to the kernel of the composition
$$\H^1(k,\Pic(A_\s)) \lra\H^1(k,\NS(A_\s))\lra\H^1(k,\Hom(A_\s,A_\s^\vee)),$$
where the second arrow is induced by the map $L\mapsto\varphi_L$.
\epr
{\em Proof.} Due to the isomorphism $\Br_{a}(A)\cong\H^1(k,\Pic(A_\s))$ 
and Proposition \ref{nuit}, the group $\Br_{a,A}(A)$ is
the kernel of the map induced by
\begin{equation}
\delta=m^*-\pi_1^*-\pi_2^*\ : \ \Pic(A_\s)\lra \Pic(A_\s\times A_\s) \label{de}
\end{equation}
on the first Galois cohomology groups $\H^1(k,-)$. By Lemma \ref{known} 
the map (\ref{de}) factors through the map of $\Ga$-modules
$\NS(A_\s)\to \Hom(A_\s,A^\vee_\s)$, where $\Hom(A_\s,A^\vee_\s)$ is a direct summand of 
the $\Ga$-module $\Pic(A_\s\times A_\s)$.  $\Box$

\medskip

M.~Stoll \cite[\S 7, p.~378]{St} denoted by $\Br_{1/2}(A)$ the subgroup of $\Br_1(A)$
consisting of the elements whose image in $\H^1(k,\Pic(A_\s))$ comes from $\H^1(k,\Pic^0(A_\s))$.
Define $$\Br_{e,1/2}(A):=\Br_e(A)\cap \Br_{1/2}(A);$$ this group is canonically isomorphic to
$$\Im\left[\H^1(k,\Pic^0(A_\s))\to \H^1(k,\Pic(A_\s))\right]
=\Ker\left[\H^1(k,\Pic(A_\s))\to \H^1(k,\NS(A_\s))\right].$$

\bco \label{coco}
Let $A$ be an abelian variety over a field $k$. 
Then $\Br_{e,1/2}(A)$ is contained in $\Br_{a,A}(A)$, and we have
$$\Br_{a,A}(A)/\Br_{e,1/2}(A)\subset \H^1(k,\NS(A_\s))[2].$$
If $\NS(A_\s)$ is a trivial $\Ga$-module, then $\Br_{a,A}(A)=\Br_{e,1/2}(A)$.
\eco
{\em Proof.} The first inclusion is immediate from Proposition \ref{rain}.
Recall that $\NS(A_\s)$ is the subgroup of $\Hom(A_\s,A_\s^\vee)$
given by the condition $\phi=\phi^\vee$. The composition of
 the inclusion $\NS(A_\s)\hookrightarrow \Hom(A_\s,A_\s^\vee)$ 
with the map $\Hom(A_\s,A_\s^\vee)\to\NS(A_\s)$
sending $\phi$ to $\phi+\phi^\vee$ is multiplication by 2.
This gives the second inclusion. 
The last statement follows since $\H^1(k,\Z)=0$ for
the trivial $\Ga$-module $\Z$. $\Box$

\bex \label{end} {\rm
Now it is easy to construct an abelian variety $A$ such that $\Br_1(A)$ is not a subgroup
of $\Br_A(A)$. Let $E$ be the elliptic curve $y^2=x^3-x$ over $k=\R$, and
let $A=E\times E$. We have direct sum decompositions of $\Ga$-modules
$$\Pic(\ov A)=\Pic(\ov E)\oplus\Pic(\ov E)\oplus \End(\ov E),$$
$$\NS(\ov A)=\NS(\ov E)\oplus\NS(\ov E)\oplus \End(\ov E),$$
and isomorphisms of $\Ga$-modules $\NS(\ov E)\cong\Z$ (with trivial $\Ga$-action) and
$\End(\ov E)\cong\Z[\sqrt{-1}]$, where 
the generator of $\Ga\cong\Z/2$ acts as complex conjugation. 
In particular, $\H^1(k,\End(\ov E))\cong\Z/2$ is a subgroup of $\H^1(k,\Pic(\ov A))$
which maps isomorphically onto $\H^1(k,\NS(\ov A))$. Next, the $\Ga$-module
$\Hom(\ov A,\ov A^\vee)\cong \End(\ov A)$ is the algebra of $(2\times 2)$-matrices over $\Z[\sqrt{-1}]$.
The Rosati involution associated to the canonical principal polarisation of $A$ acts on $\Mat_2(\Z[\sqrt{-1}])$
as the composition of transposition and complex conjugation. Hence the
injective image of $\NS(\ov A)$ consists of the matrices with diagonal entries in $\Z$ such that
the two non-diagonal entries are conjugate. The $\Ga$-submodule $\Mat_2(\Z[\sqrt{-1}])$ is the direct
sum of $\NS(\ov A)$ and the subgroup of upper-triangular matrices with purely imaginary
diagonal entries. Hence $\NS(\ov A)$ is a direct summand of 
the $\Ga$-module $\Hom(\ov A,\ov A^\vee)$. 
Thus we have an element
of order 2 in $\H^1(k,\Pic(\ov A))$ with non-zero image in $\H^1(k,\Hom(\ov A,\ov A^\vee))$
under the map of Proposition \ref{rain}, hence $\Br_1(A)$ is not contained in $\Br_A(A)$.
Note also that in this case we have $\Br_{a,A}(A)=\Br_{e,1/2}(A)$, hence
the inclusion of $\Br_{a,A}(A)/\Br_{e,1/2}(A)$ into $\H^1(k,\NS(\ov A))[2]$ is strict.}
\eex

The following statement concerns $\Br_{e,A}(A)=\Br_e(A)\cap\Br_A(A)$, which contains
$\Br_{a,A}(A)$ as a subgroup. We are grateful to the referee for pointing it out to us.

\bpr
Let $A$ be an abelian variety over a field $k$.
Then the group $$\Br_{e,A}(A)(p')/\Br_{e,1/2}(A)(p')$$ has exponent $2$.
\epr
{\em Proof.} Every element of $\Br_{e,A}(A)(p')/\Br_{e,1/2}(A)(p')$ lifts to
an element of $\Br_{e}(A)[n]$ for some
positive integer $n$ not divisible by ${\rm char}(k)$.

The natural injective map $\mu_n\to \G_m$ transforms the spectral sequence
$$\H^i(k,\H_\et^j(A_\s,\mu_n))\Rightarrow\H_\et^{i+j}(A,\mu_n).$$ 
into the spectral sequence (\ref{ss}). Thus
the exact sequences of low degree terms of these spectral sequences give rise to the following
commutative diagram with exact rows:
$$\xymatrix{0\ar[r]&\H^1(k,\H_\et^1(A_\s,\mu_n))\ar[r]\ar[d]&
\H^2_{\et, e}(A,\mu_n)\ar[r]\ar@{-{>>}}[d]&\H_\et^2(A_\s,\mu_n)\ar[d]\\
0\ar[r]&\H^1(k,\Pic(A_\s))[n]\ar[r]&\Br_e(A)[n]\ar[r]&\Br(A_\s)&}$$
Here the subscript $e$ means the kernel of the map induced by $e\colon\Spec(k)\to A$.
By the Kummer sequence for $A_\s$, we have canonial isomorphisms
$$\H_\et^1(A_\s,\mu_n)\cong \Pic^0(A_\s)[n]\cong A^\vee[n].$$
Taking Galois cohomology gives a surjective map
\begin{equation}
\H^1(k,\H_\et^1(A_\s,\mu_n))\twoheadrightarrow \H^1(k,\Pic^0(A_\s))[n].\label{sur}
\end{equation}
Since $\H_\et^1(A_\s,\mu_n)$ is a torsion group, its image in $\Pic(A_\s)$
is contained in $\Pic^0(A_\s)$, because $\NS(A_\s)$ is torsion-free.
Thus the left vertical map in the diagram factors through (\ref{sur}).

By the Kummer sequence for $A$, the middle vertical map is surjective.
It follows that the quotient of $\Br_e(A)[n]$ by $\Br_{e,1/2}(A)[n]$
is a quotient of a subgroup of $\H_\et^2(A_\s,\mu_n)$. From the canonical isomorphism
$$\H_\et^2(A_\s,\Z/n)=\wedge^2 \H_\et^1(A_\s,\Z/n)$$
we see that $[-1]$ acts trivially on $\H_\et^2(A_\s,\mu_n)$, hence also on 
$\Br_e(A)[n]/\Br_{e,1/2}(A)[n]$. 

In particular, $[-1]$ acts on $\Br_{e,A}(A)(p')/\Br_{e,1/2}(A)(p')$ trivially.
However, as follows from the last claim of Proposition \ref{nuit}, $[-1]$ acts on 
this group at $-1$. This implies our statement. $\Box$

\bigskip

\noindent Leibniz Universit\"at Hannover,
Riemann Center for Geometry and Physics,
Welfengarten 1,
30167 Hannover, Germany
\smallskip

\noindent {\tt valloni@math.uni-hannover.de}

\bigskip

\noindent Department of Mathematics,
University of Manchester, Alan Turing Building,
Oxford Road, Manchester M13 9PL England, U.K. 

\smallskip

\noindent {\tt martin.orr@manchester.ac.uk}

\bigskip

\noindent Department of Mathematics, South Kensington Campus,
Imperial College London, SW7 2BZ England, U.K. -- and --
Institute for the Information Transmission Problems,
Russian Academy of Sciences, 19 Bolshoi Karetnyi, Moscow 127994
Russia

\smallskip

\noindent {\tt a.skorobogatov@imperial.ac.uk}

\bigskip

\noindent Department of Mathematics, Pennsylvania State University, University Park, Pennsylvania 16802 USA 

\smallskip

\noindent {\tt zarhin@math.psu.edu}


{\small
\begin{thebibliography}{Mum74}
\bibitem[Ber72]{Ber72} V.G. Berkovich. The Brauer group of abelian varieties. 
{\em Funktsional. Anal. i Prilozhen.} {\bf 6:3} (1972) 10--15 (Russian).
English translation: {\em Funct. Anal. Appl.} {\bf 6:3} (1972) 180--184.

\bibitem[C18]{C18} Y. Cao. Approximation forte pour les vari\'et\'es avec 
une action d'un groupe lin\'eaire. {\em Compositio math.} {\bf 154} (2018) 773--819.

\bibitem[C20]{C} Y. Cao. Sous-groupe de Brauer invariant et obstruction de descente it\'er\'ee.
{\em Algebra Number Theory} {\bf 14} (2020) 2151--2183.

\bibitem[C]{YC} Y. Cao. Sous-groupe de Brauer invariant pour un groupe alg\'ebrique connexe quelconque. 
\texttt{arXiv:1808.06951}

\bibitem[Ch05]{Ch05} Ch.-L. Chai. Hecke orbits on Siegel modular varieties. 
{\em Geometric methods in algebra and number theory} (F. Bogomolov, Yu. Tschinkel, eds.) 
Progress in Math. {\bf 235}, pp.~71--107. Birkh\"auser, 2005.

\bibitem[CO09]{CO09} Ch.-L. Chai and F. Oort. Moduli of abelian varieties and $p$-divisible groups. 
{\em Arithmetic  Geometry}, Clay Mathematics Institute Proceedings Series {\bf 8} (H. Darmon, D. Ellwood, B. Hassett, Yu. Tschinkel, eds.), pp.~441--536. Amer. Math. Soc., 2009.

 \bibitem[COU01]{COU01} L. Clozel, H. Oh, and E. Ullmo. Hecke operators and equidistribution of Hecke points. {\em Invent. math.} {\bf 144} (2001) 327--351.

\bibitem[CTS]{CTS} J.-L. Colliot-Th\'el\`ene and A.N. Skorobogatov. {\em The Brauer--Grothendieck group.}
Springer, to appear.

\bibitem[Con06]{Con06} B. Conrad. Chow's $K/k$-image and $K/k$-trace, and the Lang--N\'eron theorem. 
{\em Enseign. Math.} {\bf 52} (2006) 37--108.

\bibitem[CR62]{CR62} C. Curtis and I. Reiner. {\em Representation theory of finite groups and associative
algebras}, 2nd ed., Interscience,  1966.

\bibitem[EGM]{vdGM} B. Edixhoven, G. van der Geer and B. Moonen. {\em Abelian varieties}.
 Preliminary version, available at https://www.math.ru.nl/$\sim$bmoonen/research.html


\bibitem[FT91]{FT} A. Fr\"ohlich and M.J. Taylor. {\em Algebraic number theory}.
Cambridge Studies in Advanced Mathematics {\bf 27}, Cambridge University Press, 1991.

\bibitem[Gro68]{Gro68} A. Grothendieck.
Le groupe de Brauer, I, II, III. {\it Dix expos\'es sur la cohomologie
des sch\'emas}. North-Holland, 1968, pp.~46--188.

\bibitem[Kem91]{Kem91} G.R.~Kempf. {\em Complex abelian varieties and theta functions}.
Universitext. Springer-Verlag, 1991.

\bibitem[Kle05]{kleiman} S.L. Kleiman. The Picard scheme. Ch.~5 of {\em Fundamental algebraic  geometry, 
Grothendieck's FGA explained}. Math. Surveys Monogr. {\bf 123}, Amer. Math. Soc., 2005.

\bibitem[Lan83]{Lang} S. Lang. {\em Complex Multiplication}.
Grundlehren der Mathematischen Wissenschaften {\bf 255}. Springer-Verlag, 1983.

\bibitem[LO98]{LO98} K.-Zh. Li and F. Oort. {\em Moduli of supersingular abelian varieties}. 
Lecture Notes in Mathematics {\bf 1680}, Springer-Verlag, 1998.

\bibitem[Mil80]{EC} J.S. Milne. {\em \'Etale cohomology}. Princeton University Press, 1980.

\bibitem[Mum74]{M} D. Mumford. {\em Abelian varieties}. 2nd ed., Oxford University Press, 1974.

\bibitem[Oort88]{O88} F. Oort. Endomorphism algebras of abelian varieties. 
{\em Algebraic geometry and commutative algebra}, Vol.~II, Konokuniya, Tokyo, 1988,
pp. 469--502.

\bibitem[OZ95]{OZ95} F. Oort and Yu.G. Zarhin. Endomorphism algebras of complex tori.
{\em Math. Ann.} {\bf 303} (1995) 11--29.

\bibitem[OSZ]{OSZ} M. Orr, A.N. Skorobogatov and Yu.G. Zarhin. 
On uniformity conjectures for abelian varieties and K3 surfaces. 
{\em Amer. J. Math.}, to appear. \texttt{arXiv:1809.01440}

\bibitem[PR91]{PR} V.P. Platonov and A.S. Rapinchuk.
{\em Algebraic groups and number theory.} Nauka, 1991. 
English translation: Academic Press, 1994.

\bibitem[Sch92]{Sch92} C. Schoen. Produkte Abelscher Variet\"aten und Moduln 
\"uber Ordnungen. {\em J. reine angew. Math.} {\bf 429} (1992) 115--123.

\bibitem[SGA3]{SGA3} {\it Sch\'emas en groupes} (SGA~3).  
S\'em. g\'eom\'etrie alg\'ebrique du Bois-Marie 1962--1964. 
Dirig\'e par M. Demazure et A. Grothendieck. 
Lecture Notes in Mathematics {\bf 151}, {\bf 152}, {\bf 153},   
Springer-Verlag, 1970. 

\bibitem[Sch05]{Sch} S. Schr\"oer. Topological methods for complex-analytic Brauer groups. 
{\em Topology} {\bf 44} (2005) 875--894.

\bibitem[Ser77]{Ser77} J.-P. Serre.
{\em Linear representations of finite groups.}
Graduate Texts in Mathematics {\bf 42}. Springer-Verlag, 1977.

\bibitem[Shi63]{Shi63} G. Shimura. On analytic families of polarized abelian varieties and automorphic functions.
{\em Ann. of Math.} {\bf 78} (1963) 149--192.

\bibitem[SZ08]{SZ08} A.N. Skorobogatov and Yu.G. Zarhin. A finiteness theorem for
the Brauer group of abelian varieties and K3 surfaces. {\em J. Alg. Geom.}
{\bf 17} (2008) 481--502.

\bibitem[SZ14]{SZ14} A.N. Skorobogatov and Yu.G. Zarhin. 
The Brauer group and the Brauer--Manin set of products of varieties. 
{\em J. Eur. Math. Soc.} {\bf 16} (2014) 749--769.

\bibitem[Sto07]{St} M. Stoll. Finite descent obstructions and rational points on curves.
{\em Algebra Number Theory} {\bf 1} (2007) 349--391.

\bibitem[Swi01]{SD} H.P.F.~Swinnerton-Dyer. {\em A brief guide to algebraic number theory}.
London Mathematical Society Student Texts {\bf 50}, Cambridge University Press, 2001.

\bibitem[Tat66]{Tat66} J. Tate. Endomorphisms of abelian varieties over finite fields.
{\em Invent. math.} {\bf 2} (1966) 134--144.


\bibitem[Zar91]{Z91} Yu.G.~Zarhin. Actions of semisimple Lie groups and orbits of Cartan subgroups.
{\em Arch. Math.} {\bf 56} (1991) 491--496.

\end{thebibliography}
}
\end{document}